\documentclass{amsart}

\usepackage{amscd}
\usepackage{amssymb}

\usepackage[T1]{fontenc}
\usepackage[latin1]{inputenc}
\usepackage{hyperref}
\usepackage[arrow,curve,matrix]{xy}
\usepackage{times}
\usepackage{graphicx}
\usepackage{color}
\usepackage{flafter}
\usepackage{placeins}
\usepackage{enumitem}

\CompileMatrices

\theoremstyle{plain}

\newtheorem{theorem}{Theorem}[section]
\newtheorem{proposition}[theorem]{Proposition}
\newtheorem{corollary}[theorem]{Corollary}
\newtheorem{lemma}[theorem]{Lemma}
\newtheorem{problem}[theorem]{Problem} 
\newtheorem{question}[theorem]{Question}
\newtheorem{conjecture}[theorem]{Conjecture}
\newtheorem{conjecture/question}[theorem]{Conjecture/Question}

\theoremstyle{definition}    
\newtheorem{definition}[theorem]{Definition}
\newtheorem{remark}[theorem]{Remark}
\newtheorem{example}[theorem]{Example}
   
\newtheorem{exercise}[theorem]{Exercise}

\theoremstyle{remark}

\def\R{\mathbf{R}}
\def\SS{\mathcal{S}}
\def\D{\mathbf{D}}

\newcommand{\PP}{\mathbf{P}}
\newcommand{\ZZ}{\mathbf{Z}}
\newcommand{\NN}{\mathbf{N}}

\newcommand{\QQ}{\mathbf{Q}}
\newcommand{\GG}{\mathbf{G}}
\newcommand{\OO}{\mathcal{O}}
\newcommand{\II}{\mathcal{I}}
   
\newcommand{\E}{\mathcal{E}}
\newcommand{\F}{\mathcal{F}}
\newcommand{\LL}{\mathcal{L}}

\newcommand{\Quot}{\mbox{Quot}_{k,d}(E)}

\newlength{\sectiontitlewidth}
\begin{document}

\title{\bf Generalized theta linear series on moduli spaces of vector bundles on curves}

\author{Mihnea Popa}
\address{Department of Mathematics \\
  University of Illinois at Chicago\\ 851 S. Morgan St., Chicago, IL 60607-7045,
  USA}
\email{mpopa@math.uic.edu}

\thanks{The author was partially supported by the NSF grant DMS-0500985, by an AMS Centennial Fellowship, and by a Sloan Fellowship, at various stages in the preparation of this work.}


\maketitle
\markboth{Mihnea Popa}
{Generalized theta linear series on moduli spaces of vector bundles on curves}

\tableofcontents

\setlength{\parskip}{.1in}

\section{Introduction}

This article is based on lecture notes prepared for the August 2006 Cologne Summer School. 
The notes have expanded in the meanwhile into a survey whose 
main focus is the study of linear series on moduli spaces of vector bundles on curves.

The paper is essentially divided into two parts. The first, smaller one,  consists of \S2 and \S3, and is almost entirely devoted to background material, reflecting precisely what I was able to cover during the 
first couple of introductory lectures at the school. There are a few short proofs or hints, some exercises, while for the fundamental results I have to refer to the standard references. This part is intended as a warm-up for the rest of the paper, and especially as an introduction to the literature for beginners in the field. The experts are advised to skip to \S4 (with perhaps the exception of \S2.3, 2.4).

The second part consists of \S4-\S7, and is a survey of the current status in the theory of linear
series and generalized theta divisors on these moduli spaces, including the connection with Quot 
schemes. Beauville's survey \cite{beauville2} has been a great source of inspiration for researchers in this field (and in fact the main one for myself). A number of problems proposed 
there have been however (at least partially) solved in recent years, and new conjectures have emerged. An update is provided in \cite{beauville5}, where the emphasis is specifically on the theta map. Here my goal is to complement this by emphasizing effective results on pluri-theta linear series, and proposing some questions 
in this direction.

Concretely, I would like to stress relatively new techniques employed in the analysis of linear 
series on moduli spaces of vector bundles, namely the use of moduli spaces 
of stable maps for understanding Quot schemes, and of the Fourier-Mukai 
functor in the study of coherent sheaves on Jacobians coming from generalized theta divisors. In addition, I will briefly describe recent important developments, namely the proof of the Strange Duality conjecture by Belkale and Marian-Oprea, and the algebro-geometric derivation of the Verlinde formula via intersection theory on Quot schemes, due to Marian-Oprea.

There is certainly recent and not so recent literature that I have not covered in this 
survey, for which I apologize. This is mainly due to the fact that in goes in directions different 
from the main thrust here, but sometimes also to my own lack of sufficient familiarity with the 
respective results. For those so inclined, much of the material can also be phrased in the language of stacks, with little or no difference in the proofs of the results in \S3-7.

\noindent
{\bf Acknowledgements.} I would like to warmly thank the organizers of the Cologne Summer School, 
Thomas Eckl and Stefan Kebekus, for an extraordinarily well-organized and successful school, which I consider a model for such events. I am also grateful to all the students who took part for their enthusiasm and their helpful comments on the material. Finally, I thank Arnaud Beauville, Prakash Belkale, Alina Marian and Drago\c s Oprea for very useful comments and corrections.

\section{Semistable bundles}

Most of the foundational material in this section can be found in detail in the books of Le Potier \cite{lepotier2} and Seshadri \cite{seshadri}, especially the construction of moduli spaces, which is only alluded to here. I have only included sketches of simpler arguments, and some exercises, in order 
to give the very beginner an idea of the subject. However, an explicit description of two special classes of vector bundles that might be of interest to the more advanced is given in \S2.3 and \S2.4.

\subsection{Arbitrary vector bundles}
Let $X$ be a smooth projective curve of genus $g$ over an algebraically closed field $k$. 
I will identify freely vector bundles with locally free sheaves.

\begin{definition}
Let $E$ be a vector bundle on $X$, of rank $r$. The \emph{determinant} of $E$ is the line bundle 
$\det E : = \wedge^r E$. The \emph{degree} of $E$ is the degree $e$ of 
$\det E$.\footnote{This is the same as the first Chern class $c_1(E)$.} The \emph{slope} of $E$ is 
$\mu(E) = \frac{e}{r} \in \QQ$.
\end{definition}

The Riemann-Roch formula for vector bundles on curves says: 
$$ \chi(E) = h^0 (X,E) - h^1 (X,E) = e + r(1-g).$$
An important example is the following: $\chi(E) = 0 \iff \mu(E) = g-1$. 

\begin{exercise}\label{minimal_degree}
Let $E$ be a vector bundle of rank $r$ on the curve $X$. Then, for each $k$,  the degree of the 
quotient bundles of $E$ of rank $k$ is bounded from below. 
\end{exercise}

From the moduli point of view, the initial idea would be to construct an algebraic variety (or scheme) parametrizing the isomorphism classes of all vector bundles with fixed invariants, i.e. rank $r$ and 
degree $e$. Note that fixing these invariants is the same as fixing the Hilbert polynomial of $E$.

\begin{definition}\label{bounded}
Let $\mathcal{B}$ be a set of isomorphism classes of vector bundles. We say that $\mathcal{B}$ is \emph{bounded}
if there exists a scheme of finite type $S$ over $k$ and a vector bundle $F$ on $S\times X$ such that all 
the elements of $\mathcal{B}$ are represented by some $F_s : = F_{|\{s\}\times X}$ with $s\in S$.
\end{definition}

\noindent
We find easily that the initial idea above is too naive.

\begin{lemma}
The set of isomorphism classes of vector bundles of rank $r$ and degree $e$ on $X$ is not bounded.
\end{lemma}
\begin{proof}
Assuming that the family is bounded, use the notation in Definition \ref{bounded}.
As $F$ is flat over $S$, by the semicontinuity theorem we get that there are only a finite number of possible $h^i F_s$ for $i = 0,1$.
Fix now a point $x\in X$ and define for each $k \in \NN$ the vector bundle
$$E_k : = \OO_X(-kx) \oplus \OO_X( (k+ e)x) \oplus \OO_X^{\oplus r-2}.$$
They all clearly have rank $r$ and degree $e$. On the other hand, when $k \to \infty$
we see that $h^0 E_k$ and $h^1 E_k$ also go to $\infty$, which gives a contradiction. 
\end{proof}

There exists however a well-known bounded moduli problem in this context, which produces the \emph{Quot scheme}. Let $E$ be a vector bundle of rank $r$ and degree $e$ on $X$, and fix integers $0\le k\le r$ and $d$. We would like to parametrize all the quotients
$$E \longrightarrow Q \longrightarrow 0$$
with $Q$ a coherent sheaf of rank $k$ and degree $d$.\footnote{Recall that the rank of  
an arbitrary coherent sheaf is its rank at a general point, while the degree is its first Chern class.} We consider the \emph{Quot functor}
$$\underline{\rm Quot}_{k,d}^E : {\rm ~Algebraic~ varieties}/k \rightarrow {\rm Sets}$$
associating to each $S$ the set of coherent quotients of $E_S : = p_X^* E$ which are flat over $S$ and have 
rank $k$ and degree $d$ over each $s\in S$. This is a contravariant functor associating to $T \overset{f}{\rightarrow} S$ the map taking a quotient $E_S \rightarrow Q$ to the quotient $E_T = (f\times id)^* E_S \rightarrow (f\times id)^* Q$.

\begin{theorem}[Grothendieck, cf. \cite{lepotier2} \S4]
There exists a projective scheme $\Quot$ of finite type over $k$, which represents the functor 
$\underline{\rm Quot}_{k,d}^E$.
\end{theorem}

This means the following: there exists a ``universal quotient'' $E_{\Quot} \rightarrow \mathcal{Q}$
on $\Quot \times X$, which induces for each variety $S$ an isomorphism 
$${\rm Hom}(S,{\rm Quot}_{k,d}(E)) \cong \underline{\rm Quot}_{k,d}^E (S)$$
given by 
$$(S\overset{f}{\rightarrow} \Quot) \to (E_S \rightarrow (f\times id)^* \mathcal{Q}).$$
The terminology is: the Quot functor (scheme)  is a \emph{fine} moduli functor (space).
I will discuss more things about Quot schemes later -- for now let's just note the 
following basic fact, which is a standard consequence of formal smoothness.

\begin{proposition}[\cite{lepotier2} \S8.2]\label{deformation}
Let $E$ be a vector bundle of rank $r$ and degree $e$, and 
$$q: ~[0\rightarrow G \rightarrow E \rightarrow F \rightarrow 0]$$
a point in $\Quot$. Then:

\noindent
(1) There is a natural isomorphism $T_q \Quot \cong {\rm Hom}(G, F) (\cong H^0 (G^\vee\otimes F))$.

\noindent
(2) If ${\rm Ext}^1 (G , F) (\cong H^1 (G^\vee \otimes F) ) = 0$, then $\Quot$ is smooth at $q$.

\noindent
(3) We have
$$h^0 (G^\vee \otimes F) \ge {\rm dim}_q \Quot \ge h^0 (G^\vee \otimes F) - h^1 (G^\vee \otimes F).$$
The last quantity is $\chi (G^\vee \otimes F) = rd - ke -k(r-k)(g-1)$ (by Riemann-Roch).
\end{proposition}

\subsection{Semistable vector bundles}

To remedy the problem explained in the previous subsection, one introduces the following notion.

\begin{definition}
Let $E$ be a vector bundle on $X$ of rank $r$ and degree $e$. It is called \emph{semistable} (respectively \emph{stable})
if for any subbundle $0\neq F \hookrightarrow E$, we have $\mu(F) \le \mu(E)$ (respectively $\mu(F) < \mu(E)$). It can be 
easily checked that in the definition we can replace subbundles with arbitrary coherent subsheaves. 
\end{definition}

\begin{exercise}\label{tensor_semistable}
Let  $E$ and $F$ be vector bundles such that $\chi (E \otimes F) = 0$. If $H^0 (E \otimes F) = 0$ 
(or equivalently, by Riemann-Roch, $H^1 (E \otimes F) = 0$), then $E$ and $F$ are semistable.
\end{exercise}

\begin{exercise}\label{stable_sequence}
Let $F$ be a stable bundle on the curve $X$. For any exact sequence 
$$0 \longrightarrow E \longrightarrow F \longrightarrow G \longrightarrow 0$$ 
we have ${\rm Hom} (G, E) = 0$.
\end{exercise}

\noindent
Here are some basic properties.

\begin{proposition}\label{simple}
If $E$ and $F$ are stable vector bundles and $\mu(E) = \mu(F)$, then every non-zero homomorphism
$\phi: E \rightarrow F$ is an isomorphism. In particular ${\rm Hom}(E, E) \cong k$ (i.e. $E$ is \emph{simple}).
\end{proposition}
\begin{proof}
Say $G = {\rm Im} (\phi)$. Then by definition we must have $\mu(E) \le \mu(G) \le\mu(F)$ and wherever we 
have equality the bundles themselves must be equal. Since $\mu(E ) = \mu(F)$, we have equality everywhere, 
which implies easily that $\phi$ must be an isomorphism. Now if $\phi \in {\rm Hom}(E,E)$, by the above we 
see that $k[\phi]$ is an finite field extension of $k$. Since this is algebraically closed, we deduce that 
$\phi = \lambda \cdot {\rm Id}$, with $\lambda \in k^*$. 
\end{proof}

\begin{exercise}
Fix a slope $\mu \in \QQ$, and let $SS(\mu)$ be the category of semistable bundles of slope $\mu$. Based on Proposition \ref{simple}, show that $SS(\mu)$ is an abelian category. 
\end{exercise}

\begin{definition}
Let $E \in SS(\mu)$. A \emph{Jordan-H\" older} filtration of $E$ is  filtration
$$ 0 = E_0 \subset E_1 \subset \ldots \subset E_p = E$$
such that each quotient $E_{i+1}/E_i$ is stable of slope $\mu$.
\end{definition}

\begin{proposition}
Jordan-H\"older filtrations exist. Any two have the same length and, upon reordering, isomorphic stable 
factors.  
\end{proposition}
\begin{proof}(\emph{Sketch})
Since the rank decreases, there is a $G \subset E$ stable of slope $\mu$. This implies that $E/G$ is semistable of slope $\mu$ and we repeat the process with $E/G$ instead of $E$. The rest is a well-known 
general algebra argument.
\end{proof}

\begin{definition}\label{definition}
(1) For any Jordan-H\"older filtration $E_{\bullet}$ of $E$, we define 
$${\rm gr}(E) : = {\rm gr}(E_{\bullet}) = \bigoplus_i E_{i+1}/E_i.$$
This is called the \emph{graded object} associated to $E$ (well-defined by the above).

\noindent
(2) A vector bundle $E$ is called \emph{polystable} if it is a direct sum of stable bundles of the same 
slope. (So for $E$ semistable, ${\rm gr}(E)$ is polystable.)

\noindent
(3) Two bundles $E,F \in SS(\mu)$ are called \emph{$S$-equivalent} if ${\rm gr}(E)\cong {\rm gr}(F)$.
\end{definition}

Sometimes we can reduce the study of arbitrary bundles to that of semistable ones via the following:

\begin{exercise}[Harder-Narasimhan filtration]
Let $E$ be a vector bundle on $X$. Then there exists an increasing filtration
$$ 0 = E_0 \subset E_1 \subset \ldots \subset E_p = E$$
such that 
\begin{enumerate}
\item Each quotient $E_{i+1}/E_i$ is semistable.
\item We have $\mu (E_i/E_{i-1}) > \mu(E_{i+1}/E_i)$ for all $i$.
\end{enumerate}
The filtration is unique; it is called the \emph{Harder-Narasimhan filtration} of $E$.
\end{exercise}

\noindent
\begin{example}
(1) All line bundles are stable. Any extension of vector bundles in $SS(\mu)$ is also in $SS(\mu)$.

\noindent
(2) If $(r,e) = 1$, then stable is equivalent to semistable.

\noindent
(3) If $X = \PP^1$, by Grothendieck's theorem we know that every vector bundle splits as 
$E \cong \OO(a_1) \oplus \ldots \oplus \OO(a_r)$, so it is semistable iff all $a_i$ are equal.

\noindent
(4) We will encounter some very interesting examples below. Until then, here is the first type of example which requires a tiny argument. Say $L_1$ and $L_2$ are line bundles on $X$, with ${\rm deg}~L_1 = d$ and ${\rm deg}~L_2 = d+1$. Consider extensions of the form
$$0 \longrightarrow L_1 \longrightarrow E \longrightarrow L_2\longrightarrow 0.$$ 
These are parametrized by ${\rm Ext}^1 (L_2, L_1) \cong H^1 (X, L_1 \otimes L_2^{-1})$. By Riemann-Roch this is isomorphic to $k^{g-2}$, so as soon as $g \ge 3$ we can choose the extension to be non-split. For such a choice $E$ is stable: first note that $\mu(E) = d + 1/2$.
Consider any line subbundle $M$ of $E$. If ${\rm deg}~M\le d$ everything is fine. If not, the induced map $M \rightarrow L_2$ must be non-zero (otherwise $M$ would factor through $L_1$, of too low degree). This immediately implies that it must be an isomorphism, which is a contradiction since the extension is non-split.
\end{example}

\noindent
For later reference, let me also mention the following important result. An elementary proof based on Gieseker's ideas can be found in \cite{positivity} \S6.4.

\begin{theorem}\label{tensor}
Assume that ${\rm char}(k) =0$. If $E$ and $F$ are semistable bundles, then $E\otimes F$ is also 
semistable.\footnote{It is also true that if $E$ and $F$ are actually stable, then $E\otimes F$ is 
polystable.} In particular, for any $k$, $S^k E$ and $\wedge^k E$ are also semistable.
\end{theorem}

\subsection{Example: Lazarsfeld's bundles}
Here are the more interesting examples of semistable bundles promised above. In this subsection, following \cite{lazarsfeld},  I describe a construction considered by Lazarsfeld in the study of syzygies of curves. In the next I will describe a different construction due to Raynaud.

Consider a line bundle $L$ on $X$ of degree $d\geq 2g+1$. 
Denote by $M_{L}$ the kernel of the evaluation map:

$$0 \longrightarrow M_{L} \longrightarrow H^{0}(L)\otimes \mathcal{O}_{X}
\overset{{\rm ev}}{\longrightarrow} L \longrightarrow 0$$
and let $Q_{L}= M_{L}^{\vee}$. Note that ${\rm rk}~Q_L = h^0L -1 = d-g$ and ${\rm deg}~ L = d$, 
so $\mu(Q_L) = \frac{d}{d-g}$. The main property of $Q_L$ is the following:

\begin{proposition}\label{filtration}
If $x_{1},\ldots,x_{d}$ are the points of a generic hyperplane section of
$X\subset {\mathbb P}(H^{0}(L))$, then $Q_{L}$ sits in an extension:
$$0 \longrightarrow 
\bigoplus_{i=1}^{d-g-1}\mathcal{O}_{X}(x_{i})
\longrightarrow
Q_{L} \longrightarrow \mathcal{O}_{X}(x_{d-g}+\ldots +x_{d})
\longrightarrow 0.$$
\end{proposition}
\begin{proof}
(1) Look first at a general situation: assume that $x_1, \ldots, x_k$ are distinct points 
on $X$, with $D = x_1 + \ldots + x_k$, such that:
\begin{enumerate}
\item $L(-D)$ is globally generated.
\item $h^0 L(-D) = h^0 L -k$.
\end{enumerate}
\emph{Claim:} In this case we have an exact sequence
$$0\longrightarrow  M_{L(-D)} \longrightarrow M_L \longrightarrow \bigoplus_{i=1}^{k}\mathcal{O}_{X}(-x_{i}) \longrightarrow 0.$$

\begin{proof}
By Riemann-Roch and (1) we have $h^1 L(-D) = h^1 L$. This implies that after passing to cohomology, 
the exact sequence 
$$0 \longrightarrow L(-D) \longrightarrow L \longrightarrow L_D \longrightarrow 0$$
stays exact. Note also that the evaluation map for $L_D$ can be written as
$$\big(H^0 L_D \otimes \OO_X \rightarrow L_D \big) = \big( \bigoplus_{i=1}^k (\OO_X\rightarrow \OO_{x_i})\big).$$
Using all of these facts, (2) and the Snake Lemma, we obtain a diagram whose top row gives precisely our Claim.

$$\xymatrix{
& 0 \ar[d] & 0 \ar[d] & 0 \ar[d] & \\
0 \ar[r] & M_{L(-D)} \ar[r] \ar[d] & M_L  \ar[r] \ar[d] & \oplus_{i=1}^{k} \OO_X (-x_{i}) \ar[r] \ar[d] & 0 \\
0 \ar[r] & H^0 L(-D) \otimes \OO_X \ar[r] \ar[d]  & H^0 L \otimes \OO_X  \ar[r] \ar [d] &  H^0 L_D \ar[r] \ar[d] & 0 \\
0 \ar[r] & L(-D) \ar[r] \ar[d]  & L  \ar[r] \ar[d] & L_D \ar[r] \ar[d] & 0 \\
& 0  & 0  & 0 & } $$ 
\end{proof}

\noindent
(2) Note now that if in the situation above we also have $h^0 L(-D) = 2$, then
the exact sequence takes the form
$$0\longrightarrow  \OO_X(x_1 + \ldots + x_k) \otimes L^{-1} \longrightarrow M_L \longrightarrow \bigoplus_{i=1}^{k}\mathcal{O}_{X}(-x_{i}) \longrightarrow 0.$$ 
This is clear, since in this case by definition ${\rm rk} ~M_{L(-D)} = 1$ and ${\rm det}~M_{L(-D)} = L^{-1} (D)$.

\noindent 
(3) Say now that $X$ is embedded in $\PP^n$ by a line bundle $L$, and $\Lambda : = {\rm Span} \big(
x_1, \ldots, x_k\big) \subset \PP^n$. Then we have the following general 

\begin{exercise}\label{linear_conditions}
(a) $h^0 L(-D) = h^0 L - k \iff {\rm dim}~\Lambda = k-1$.

\noindent
(b) $L(-D)$ is globally generated $\iff \Lambda \cap X = \{x_1, \ldots, x_k\}$
and $\Lambda$ does not contain the tangent line to $X$ at $x_i$ for all $i$.
\end{exercise}

\noindent
(4) Now consider the situation in the Proposition, when $X \subset \PP (H^0 L)$ and $x_1, \ldots, x_d$ are 
the points of a general hyperplane section of $X$. This implies that the points are in linear general 
position, which immediately gives that the conditions in Exercise \ref{linear_conditions} are satisfied. By the previous points, we are done.
\end{proof}

Proposition \ref{filtration} implies the stability of $Q_L$. The proof below is due to 
Ein-Lazarsfeld \cite{el1}.

\begin{proposition}
Under the assumptions above $Q_L$ is a stable bundle.
\end{proposition}
\begin{proof}
Let's see that the dual $M_L$ is stable.
One can actually prove a bit more: $M_L$ is \emph{cohomologically stable}, i.e. for any 
line bundle $A$ of degree $a$ and any $t < {\rm rk} M_L = d-g$: 
\begin{equation}\label{cs}
H^0 (\bigwedge^t M_L \otimes A^{-1}) = 0 {\rm ~~~if~} a \ge t\cdot \mu(M_L) = -\frac{td}{d-g}.
\end{equation}
This implies the stability of $M_L$: indeed, if $F \hookrightarrow M_L$ is a subbundle of degree 
$a$ and rank $t$, then we have an inclusion $A :=\bigwedge^t F \hookrightarrow \bigwedge^t M_L$, which
implies that $H^0 (\bigwedge^t M_L \otimes A^{-1}) \neq 0$. By cohomological stability we must
have $\mu(F) = \frac{a}{t} < \mu(M_L)$, so $M_L$ is stable.

To prove (\ref{cs}), take exterior powers in the dual of the sequence in Proposition \ref{filtration}
to obtain
$$0 \rightarrow \OO_X(-x_{d-g} - \ldots - x_d) \otimes \bigwedge^{t-1} (\bigoplus_{i=1}^{d-g-1}\OO_{X}(-x_{i})) \rightarrow \bigwedge^t M_L
\rightarrow \bigwedge^t (\bigoplus_{i=1}^{d-g-1}\OO_{X}(-x_{i})) \rightarrow 0.$$
In other words we have an exact sequence
$$0 \rightarrow \underset{1\le i_1 < \ldots < i_{t-1}\le d-g-1}{\bigoplus}\OO_{X}(-x_{i_1}-\ldots - x_{i_{t-1}} - x_{d-g}\ldots -x_d)) \rightarrow \bigwedge^t M_L $$
$$
\rightarrow  \underset{1\le j_1 < \ldots < j_{t}\le d-g-1}{\bigoplus}\OO_{X}(-x_{j_1}-\ldots - x_{j_{t}}))\rightarrow 0.$$
We tensor this sequence by $A^{-1}$. It can be checked easily that on both extremes $H^0$ is zero, 
as the points $x_i$ are general. This implies what we want.
\end{proof}

\begin{corollary}
If ${\rm char}(k)=0$, then for all $p$ the bundle $\bigwedge^p Q_L$ is semistable.
\end{corollary}

\subsection{Example: Raynaud's bundles}(cf. \cite{raynaud})
Let $X\hookrightarrow J(X)$ be an Abel-Jacobi embedding. Denote $A = J(X)$. Then $A$ has a principal polarization $\Theta$, which induces an isomorphism $A \cong \widehat A = {\rm Pic}^0(A)$.
Denote by $\mathcal{P}$ a Poincar\'e bundle on $A \times \widehat A$. For any $m \ge 1$, we consider what is called the \emph{Fourier-Mukai} transform of $\OO_{\widehat A}(-m\Theta)$, namely: 
$$F : = \widehat{\OO_{\widehat A} (-m\Theta)} : = R^g{p_A}_* (p_{\widehat A}^* \OO(-m\Theta) \otimes \mathcal{P}).$$ 
By base change, this is a vector bundle with fiber over $x \in A$ isomorphic to 
$$H^g ( \widehat A, \OO_{\widehat A} (-m\Theta) \otimes \mathcal{P}_{|\{x\}\times \widehat A}) \cong 
H^0 (J(X), \OO_{J(X)} (m\Theta) \otimes P_x)^\vee,$$
where $P_x$ is the line bundle in ${\rm Pic}^0(J(X)))$ corresponding to $x \in J(X)$. (Note that as $x$ varies with $J(X)$, $P_x$ varies with ${\rm Pic}^0(J(X))$.) Hence $F$ is a vector bundle of rank 
$m^g$. Define the vector bundle $E := F_{|X}$ on $X$.

The claim is that this is a semistable bundle. Indeed, consider the multiplication by $m$ map 
$\phi_m: A \rightarrow A$. By a result of Mukai, \cite{mukai2} Proposition 3.11, we have that 
$$\phi_m ^* F \cong H^0 \OO_{\widehat A} (m \Theta) \otimes \OO_{\widehat A} (m \Theta) 
\cong \underset{m^g}{\bigoplus} \OO_{\widehat A} (m \Theta) .$$
We consider the \' etale base change $\psi : Y \rightarrow X$, where $Y = \phi_m^{-1} (X)$.  
The decomposition of $F$ via pull-back by $\phi_m$ implies that $\psi^* E $ is semistable on $Y$. 
Applying the well-known Lemma below, we deduce that $E$ is semistable. 

\begin{lemma}(cf. \cite{positivity} Lemma 6.4.12)
Let $f: Y \rightarrow X$ be a finite morphism of smooth projective curves, and $E$ a vector bundle on $X$. Then $E$ is semistable if and only if $f^*E$ is semistable.
\end{lemma}

Let's also compute the slope of $E$. Note that ${\rm deg}~\psi^* E = {\rm deg}~\phi_m \cdot {\rm deg}~E = m^{2g} \cdot {\rm deg}~E$. This gives 
$${\rm deg}~E = \frac {m^g \cdot m \cdot (\theta\cdot [Y])}{m^{2g}} = \frac{\theta\cdot [Y]}{m^{g-1}}.$$
But since $Y = \phi_m^{-1}(X)$, $\phi_m^* \theta \equiv m^2 \cdot \theta$ and $\theta \cdot [X]= g$, we have that $m^2 \cdot (\theta \cdot [Y]) = g \cdot m^{2g}$. Putting everything together, we finally obtain 
$${\rm deg}~E = g\cdot m^{g-1}, {\rm~i.e.~} \mu(E) = \frac{g}{m}.$$

\begin{remark}
Schneider \cite{schneider2} uses a slightly more refined study of Raynaud's bundles based on theta-group actions in order to show that they are actually \emph{stable}.
\end{remark}

\begin{remark}
A generalization of Raynaud's examples, by means of semi-homogeneous vector bundles on the 
Jacobian of $X$, is presented in Remark \ref{verlinde_other} (2).
\end{remark}

\begin{remark}
I will only mention in passing one more interesting idea, which has to do with moduli of curves: certain types of stable bundles can exist only on special curves, just like in the usual Brill-Noether theory for line bundles. One can look 
for vector bundles with ``many" sections, i.e. wonder whether for a given $k$ there exist semistable bundles $E$ on $X$ of rank $r$ and degree $d$ (or fixed determinant $L$), which have $k$ independent sections. 
For instance, if $g(X) = 10$, one can see that the condition that there exist a semistable bundle 
of rank $2$ and determinant $\omega_X$ with at least $7$ sections is codimension $1$ in 
$\mathcal{M}_{10}$, i.e. such bundles exist only on curves filling up a divisor in $\mathcal{M}_{10}$. The closure of this locus in $\overline{\mathcal{M}}_{10}$ is a very interesting divisor, the first one shown to be of \emph{slope} smaller than expected (cf. \cite{fp}). 
As examples of stable bundles though, the vector bundles here are not quite new: they are roughly speaking of the same kind as the Lazarsfeld examples above.
\end{remark}

\subsection{The moduli space}
Back to the boundedness problem: we want to see that semistable bundles do the job. First a technical point, extending a well-known fact for line bundles.

\begin{exercise}\label{vanishing}
Let $E$ be a semistable bundle on $X$.

\noindent
(a) If $\mu(E) > 2g-2$, then $h^1 E =0$.

\noindent
(b) If $\mu(E) > 2g-1$, then $E$ is globally generated. 
\end{exercise}

\begin{proposition}
The set $\mathcal{S}(r,d)$ of isomorphism classes of semistable bundles of rank $r$ and 
degree $d$ is bounded.
\end{proposition}
\begin{proof}
Fix $\OO_X(1)$ a polarization on the curve. By Exercise \ref{vanishing}, there exists a fixed $m \gg 0$ such that for all $F$ in $\mathcal{S}(r,d)$ we have $h^1 F(m) = 0$ and $F(m)$ is globally generated.
Let $q := h^0 F(m) = \chi (F(m))$, which is constant by Riemann-Roch. The global generation of $F(m)$
means that we have a quotient
$$\OO_X^{\oplus q} (-m) \overset{\beta}{\longrightarrow} F \longrightarrow 0.$$
These all belong to the Quot scheme ${\rm Quot}_{r,d}(\OO_X^{\oplus q} (-m))$, which is 
a bounded family.
\end{proof}

The quotient $\beta$ can be realized in many ways: fix a vector space $V \cong k^q$, and choose an
isomorphism $V \cong H^0 F(m)$. On ${\rm Quot}_{r,d} (V \otimes \OO_X(-m)) \cong 
{\rm Quot}_{r,d} (\OO_X^{\oplus q}(-m))$ we have a natural $GL(V)$-action, namely each $g \in GL(V)$ 
induces a diagram
$$\xymatrix{
V \otimes \OO_X(-m) \ar[r] \ar[d]_{g\otimes {\rm id}} & Q  \\
V \otimes \OO_X(-m) \ar[ur] &  } $$ 
The scalar matrices act trivially on these quotients, so in fact we have a $PGL(V)$-action.

\begin{proposition}
Let $\Omega \subset {\rm Quot}_{r,d} (V \otimes \OO_X(-m))$ be the set of quotients $Q$ such that
$Q$ is semistable and $V \cong H^0 Q(m)$. Then $\Omega$ is invariant under the $PGL(V)$-action, and 
we have a bijection
$$\Omega / PGL(V) \longrightarrow \mathcal{S}(r,d).$$ 
\end{proposition}
\begin{proof}
The set $\Omega$ is clearly invariant, and the points in the same orbit give isomorphic quotient 
bundles, so we have a natural map $\Omega / PGL(V) \longrightarrow \mathcal{S}(r,d)$. Since $m \gg 0$, 
by global generation the map is surjective. Suppose now that we have two different quotients inducing
an isomorphism $\phi$:
$$\xymatrix{
V \otimes \OO_X(-m) \ar[r]  & Q  \ar[d]_{\phi}^{\cong}\\
V \otimes \OO_X(-m) \ar[r] &  Q} $$ 
This induces an isomorphism $H^0 Q(m) \overset{H^0 \phi(m)}{\longrightarrow} H^0 Q(m)$, which 
corresponds to an element $g \in GL(V)$. This is uniquely determined up to scalars, so we can in fact 
consider it in $PGL(V)$.
\end{proof}

Although this was not one of the topics of the lectures, the main point of Geometric 
Invariant Theory (GIT) -- appropriately for this volume, one of D. Mumford's most celebrated 
achievements --  is  in this context essentially to show that in fact this quotient has the structure 
of a projective algebraic variety.\footnote{This is literally true if we consider $S$-equivalence classes 
of semistable bundles (cf. Definition \ref{definition}) instead of isomorphism classes.} The GIT machinery constructs the space $U_X(r,d)$ -- an algebraic variety replacement of $\mathcal{S}_{r,d}$ --  
the \emph{moduli space of $S$-equivalence classes of semistable vector bundles of rank $r$ and degree $d$ on $X$}. For the well-known construction, that can be adapted to arbitrary characteristic 
without difficulty, cf. \cite{lepotier2} \S7 or \cite{seshadri} Ch.I.
We denote by $U_X^s(r,d)$ the open subset corresponding to isomorphism classes of stable bundles.  

We consider also a variant of $U_X(r,d)$ when the determinant of the vector bundles is fixed. More precisely, for any $L \in {\rm Pic}^d(X)$, we denote by $SU_X(r, L)$ the moduli space of ($S$-equivalence 
classes of) semistable bundles of rank $r$ and fixed determinant $L$. These are the fibers of the 
natural determinant map
$${\rm det}: U_X(r,d) \longrightarrow {\rm Pic}^d(X).$$
These spaces have the following basic properties, essentially all derived from the GIT construction 
(cf. e.g. \cite{seshadri} Ch.I and \cite{dn}):

\noindent
(1) $U_X(r,d)$ and $SU_X(r, L)$ are projective varieties (i.e irreducible), of dimension $r^2 (g-1) + 1$
and $(r^2 -1)(g-1)$ respectively.

\noindent
(2) $U_X(r,d)$ and $SU_X(r,L)$ are in general only \emph{coarse} moduli spaces (i.e. there is no universal family). In fact one can show that they are fine moduli spaces if and only if $(r,d) =1$.

\noindent
(3) $U_X(r,d)$ and $SU_X(r, L)$ are normal, Gorenstein, with rational singularities.

\noindent
(4) ${\rm Sing} (U_X(r,d)) = U_X(r,d) - U_X^s(r,d)$, unless $g = 2$, $r=2$ and $d =$ even, when $U_X(r,d)$ is smooth. Same for $SU_X(r, L)$. 

\noindent
(5) If $E$ is stable, then $T_E U_X(r,d) \cong H^1 (X, E^\vee \otimes E) \cong {\rm Ext}^1(E, E)$ and 
$T_E SU_X(r, L) \cong H^1 (X, {\rm ad}(E))$. (Recall that ${\rm ad}(E)$ is the sheaf of trace $0$ 
endomorphisms of $E$.)

\begin{exercise}
Let $r > 0$ and $d$ be two integers, $h := (r,d)$ and $r_0 : = r/h$, $d_0 := d/h$. The dimension of the
singular locus ${\rm Sing}~U_X(r,d)$ is equal to $\frac{r^2}{2} (g-1) + 2$  if $h$ is even, and 
$\frac{r^2 + r_0^2}{2} (g-1) + 2$ if $h$ is odd.
\end{exercise}

\begin{exercise}[Moduli of vector bundles on elliptic curves] 
Let $X$ be an elliptic curve. Assume that $(r,d) = 1$ and $d > 0$. Then the determinant map ${\rm det}: U_X(r,d) \rightarrow {\rm Pic}^d(X)$ is an isomorphism. 
(Cf. \cite{lepotier2}, Exercise 8.7, for a break-down of this exercise. The case when $(r,d) >1$ is more involved. It can be shown that there are no stable bundles of rank 
$r$ and degree $d$, and that $U_X(r,d)$ is isomorphic to the symmetric product $S^h X$, where 
$h = (r,d)$; cf. \cite{atiyah}.)
\end{exercise}

\section{Generalized theta divisors}

Consider first the following well-known situation: on ${\rm Pic}^d (X)$ we have, for each $L \in {\rm Pic}^{g-1-d}(X)$,  a theta divisor
$$\Theta_L : = \{ M~|~ h^0(M\otimes L) \neq 0\},$$
which is a ``translate'' of the principal polarization $\Theta$ on the Jacobian $J(X) \cong {\rm Pic}^0(X)$
(or of $W_{g-1} = \{N ~|~h^0 N \neq 0\} \subset {\rm Pic}^{g-1}(X)$).
Note that the numerical choice is such that $\chi (M \otimes L) = 0$.

Now fix $E \in U_X(r,d)$.\footnote{I will always  use somewhat abusively vector bundle notation instead of $S$-equivalence class notation for simplicity. For anything we are interested in, it can be checked easily using Jordan-H\"older filtrations that the statement is independent of the choice in the $S$-equivalence class.} When can we have a vector bundle $F$ such that $\chi(E \otimes F) = 0$? By Riemann-Roch we need $\mu(E \otimes F) = g-1$, in other words
$$\mu(F) = g-1 - \mu(E).$$
Using the notation  $h = (r,d)$, $r_0 = r/h$ and $d_0 = d/h$, we see that the only possibilities are 
${\rm rk}~F = kr_0$ and ${\rm deg}~F = k(r_0(g-1) -d_0)$, with $k \ge 1$. Fix such an $F$. 

\noindent
\emph{Claim:} If there exists $E \in U_X(r,d)$ such that $H^0 (E \otimes F) = 0$, then 
$$\Theta_F : = \{ E ~|~ h^0 (E \otimes F ) \neq 0\} \subset U_X(r,d)$$
is a divisor with a natural scheme structure (a \emph{generalized theta divisor}). The same is true on 
$SU_X(r, L)$.

\begin{proof}
I will sketch here only the case $(r,d)= 1$, when there exists a universal family on $X \times U_X(r,d)$.\footnote{The general case follows by doing essentially the same thing on a Quot scheme and using GIT, or on an \'etale cover of $U_X(r,d)$ (cf. \cite{dn}).}
It is a typical example of a determinantal construction: let $\E$ be a universal bundle on $X \times U_X(r,d)$ and $D$ an effective divisor on $X$ with 
${\rm deg}~D \gg 0$.  We can consider the following natural sequence obtained by pushing forward to 
$U_X(r,d)$: 
$$0\rightarrow  {p_U}_* (\E \otimes p_X^* F)\rightarrow  {p_U}_* (\E \otimes p_X^* F(D))\rightarrow  {p_U}_* ((\E \otimes p_X^* F(D))_{|D \times U_X(r,d)})\rightarrow $$
$$\rightarrow  R^1{p_U}_* (\E \otimes p_X^* F)\rightarrow 0,$$
Note that the $0$ on the right is obtained by Base Change, since for any $E \in U_X(r,d)$ we have 
that $h^1(E \otimes F (D)) = 0$, as $D$ has sufficiently large degree and the family of $E$'s is bounded. Let's redenote this sequence 
$$0 \rightarrow K \rightarrow G \overset{\alpha}{\rightarrow} H \rightarrow C \rightarrow 0.$$
In fact we have $K = 0$ (exercise: why?). In any case, by Base Change and Riemann-Roch, $G$ and $H$ are vector bundles on $U_X(r,d)$ of the same rank $r \cdot {\rm rk}(F) \cdot {\rm deg}(D)$. On the locus where Base Change applies, fiberwise we have 
$$H^0 (E \otimes F(D) ) \rightarrow H^0 (E \otimes F(D)_D) \rightarrow 
H^1 (E \otimes F).$$
We see immediately that the degeneracy locus of $\alpha$ is set-theoretically precisely $\Theta_F$. 
This means that $\Theta_F$ has a determinantal scheme structure and since ${\rm rk} ~G = 
{\rm rk}~H$ we have that ${\rm codim}~\Theta_F \le 1$. Since there exists $E$ such that 
$h^0 (E\otimes F) = 0$, i.e $E \not\in \Theta_F$, we must then have ${\rm codim}~\Theta_F = 1$.
\end{proof}

Note that by Exercise \ref{tensor_semistable}, a vector bundle $F$ as in the Claim has to be semistable. The requirement above is indeed satisfied for $F$ general in its corresponding moduli space of semistable bundles, so we always have enough generalized theta divisors (cf. e.g. \cite{dn} 0.2). One proof is based on counting techniques similar to those explained below in \S7.
Note however that, by semicontinuity, it is obvious for example on the moduli spaces that contain direct sums of  line bundles, for instance $U_X(k, k(g-1))$.

One of the most important facts that has been proved about the moduli spaces of vector bundles on curves, again using the GIT description, is the following: 

\begin{theorem}(Drezet-Narasimhan \cite{dn})\label{picard}
(1) $U_X(r,d)$ and $SU_X (r, L) $ are locally factorial.

\noindent
(2) For any $F \in  U_X(k r_0, k(r_0(g-1) -d_0))$ such that $\Theta_F$ is a divisor, the line bundle
$\OO(\Theta_F)$ on $SU_{X}(r,L)$ does not depend on the choice  of $F$.
The Picard group of $SU_{X}(r,L)$ is isomorphic to $\ZZ$, generated by an ample line bundle 
$\LL$ (called the \emph{determinant line bundle}), and $\Theta_F \in |\LL^k|$.

\noindent
(3) The inclusions ${\rm Pic}({\rm Pic}^d(X))\subset {\rm Pic}(U_{X}(r,d))$ (given by the determinant morphism ${\rm det}: U_X(r,d) \rightarrow {\rm Pic}^d (X)$) and $\ZZ\cdot\OO(\Theta_F)\subset {\rm Pic}(U_{X}(r,d))$, with $k =1$, induce an isomorphism
$${\rm Pic}(U_{X}(r,d))\cong {\rm Pic}({\rm Pic}^d(X))\oplus\ZZ.$$
We have the following transformation formula: if $F, F^\prime \in U_X(k r_0, k(r_0(g-1) -d_0))$, then
$$\OO(\Theta_F) \cong \OO(\Theta_{F^\prime}) \otimes {\rm det}^* ( {\rm det}(F) \otimes {\rm det}(F^\prime)^{-1}).$$

\noindent
(4) $\omega_{SU_X(r, L)} \cong \LL^{-2h}$, where $h = (r,d)$ (in particular, $SU_X(r, L)$ are 
Fano varieties). There is a similar formula for $U_X(r,d)$.
\end{theorem}

\begin{exercise}\label{pullback}
For $L \in {\rm Pic}^d (X)$, consider the tensor product map 
$$\tau: SU_{X}(r, L)\times {\rm Pic}^0(X)  \longrightarrow U_{X}(r,d).$$
\noindent
(1) Show that $\tau$ is \'etale Galois, with Galois group $X[r] \subset {\rm Pic}^0 (X)$, the group of 
$r$-torsion points in the Jacobian of $X$. 

\noindent
(2) For any vector bundle $F \in U_X(kr_0, k(r_0(g-1) - d_0))$ giving a theta divisor $\Theta_F$ on $U_X(r,d)$, show that 
$$\tau^* \OO_U (\Theta_F) \cong \LL^k \boxtimes \OO_{{\rm Pic}} (krr_0\Theta_N),$$
where $\LL$ is the determinant line bundle on $SU_X(r,L)$ and $N \in {\rm Pic}^{g-1}(X)$.
\end{exercise}

The very existence of generalized theta divisors has as a consequence a useful
criterion for basepoint-freeness, point separation and tangent vector separation. 
I state it for simplicity only in degree $0$, but there are analogous statements for any degree.
The first two parts follow immediately from the definition of $\Theta_F$.

\begin{corollary}\label{effective}
(1) Given $E \in SU_X(r)$, if there exists $F \in U_X(k, k(g-1))$ such that $H^0 (E \otimes F ) = 0$, 
then $|\LL^k|$ is free at $E$.

\noindent
(2) Given $E_1, E_2 \in SU_X(r)$, if there exists $F \in U_X(k, k(g-1))$ such that $H^0 (E_1 \otimes F ) = 0$ and $H^0 (E_2 \otimes F ) \neq 0$, 
then $|\LL^k|$ separates $E_1$ and $E_2$.

\noindent
(3) Consider $E \in SU_X^s(r)$ and $v \in T_E SU_X (r)$, corresponding to an extension
$$0 \longrightarrow E \longrightarrow \E \longrightarrow E \longrightarrow 0.$$
If there exists $F \in U_X(k, k(g-1))$ such that $h^0 (E \otimes F ) = 1$ and the induced element 
in ${\rm Hom}(F^\vee, E)$ does not lift to an element in ${\rm Hom}(F^\vee, \E)$ via the extension map, 
then $|\LL^k|$ separates $v$.
\end{corollary}

\begin{exercise}
Prove part (3) of Corollary \ref{effective}.
\end{exercise}

\begin{remark}
Among the generalized theta divisors described in this section, an important class is provided by pluri-theta divisors 
on Jacobians associated to higher rank vector bundles. This means, say for $E$ in $U_X(r, r(g-1))$, loci of the form
$$\Theta_E = \{ L \in {\rm Pic}^0 (X) ~|~ H^0 (E \otimes L) = 0\} \subset {\rm Pic}^0 (X).$$
For a variety of reasons, it is interesting to detect when a bundle $E$ ``has" such a theta divisor, i.e. when $\Theta_E$ is indeed a proper subset. Besides the connection with linear series on the moduli space described in these notes, some special such theta divisors appear for example in \cite{fmp} in connection with the Minimal Resolution Conjecture for points on curves. Cf. also \cite{beauville3} for related topics.
\end{remark}

\section{Quot schemes and stable maps}

The Quot scheme is a key tool in the study of moduli spaces, beyond its obvious role 
coming from the GIT construction. For applications, one of the most important 
things is to understand a good dimension bound. If $X$ is a smooth projective curve, and $E$ is a vector bundle on $X$, denote by $\Quot$ the Quot scheme parametrizing coherent quotients of $E$ of rank $k$ and degree $d$. Set 
$$d_k = d_k (E):=\underset{{\rm rk}\,F=k}{{\rm min}}\left\{{\rm deg}\,F~|\,\,\,\mbox{$F$ 
a quotient of $E$}\right\},$$  
i.e. the minimal degree of a quotient bundle of rank $k$ (cf. Exercise \ref{minimal_degree}).

\begin{theorem}[\cite{pr}]\label{dimension}
For an arbitrary vector bundle $E$ of rank $r$, the following holds: 
$$\dim {\rm Quot}_{k,d}(E) \leq k(r-k)+ (d-d_{k})r, ~for~ all~ d\ge d_k.$$
\end{theorem}

This result opens the door towards obtaining effective bounds for basepoint-freeness and very ampleness on the moduli spaces of vector bundles, as explained in \S7, and is useful for a variety of 
counting problems related to families of vector bundles. 

\begin{remark}\label{subbundles}
(1) In general the invariant $d_k$ is quite mysterious. In case $E$ is semistable however, there are obvious lower bounds given by the semistability condition: for example if $E\in SU_X(r)$, then 
$d_k \ge 0$.

\noindent
(2) Lange \cite{lange} showed that if $E\in U_X(r,e)$ is general, then $d_k (E)$ is the smallest integer making the
expression $d_k r - ke - k(r-k) (g-1)$ non-negative.

\noindent
(3) By Proposition \ref{deformation}, if $q = [E \rightarrow Q \rightarrow 0]$ is a quotient of rank $r$ and degree $d$, we always have 
$${\rm dim}_q {\rm Quot}_{k,d}(E) \ge rd - ke - k(r-k)(g-1).$$
Combining with the Theorem above, this gives 
$$rd_k - ke \le k(r-k)g,$$
which is an upper bound on $d_k$ due to Lange \cite{lange} and Mukai-Sakai \cite{ms}, 
extending an old theorem of Nagata on the minimal self-intersection of a section of a ruled surface.
Loosely speaking, every vector bundle has quotients of every rank of sufficiently small degree.
\end{remark} 

I will explain briefly the approach to proving Theorem \ref{dimension}. 
The starting point is the well-known identification provided by the following:

\begin{exercise}
Let $E$ be a vector bundle of rank $r$ and degree $e$ on a smooth projective curve $X$. 
  
\noindent
(1) Show that all the vector bundle quotients of $E$ of rank $k$ can be identified with sections 
of the Grassmann bundle (of quotients) projection $\pi : \GG (k , E) \rightarrow X$.   

\noindent
(2) Let $E \overset{q}{\rightarrow} F \rightarrow 0$ be a vector bundle quotient of $E$ of rank $k$. 
By the above this corresponds to a section $s: X \rightarrow \GG (k , E)$. If $Y = s(X)$, then 
show that 
$${\rm deg}~F = d \iff \OO_{\GG}(1)\cdot Y = rd -ke$$ 
where $\OO_{\GG}(1)$ on $\GG (k , E)$ is the line bundle inducing the Pl\"ucker line bundle on the fibers. Deduce that $d$ determines and is determined by the cohomology class of $Y$.
\end{exercise}

It had already been noticed by Mukai-Sakai \cite{ms} that the Rigidity 
Lemma applied to sections as in the Exercise above provides the bound in the case $d = d_k$. The general procedure is an induction relying on the Bend-and-Break Lemma, which needs to be applied to a close kin of $\Quot$, namely the moduli space of stable maps
$\overline{\mathcal{M}}_g (\GG (k, E), \beta_d)$. Here $\beta_d$ is the class of any section corresponding to a quotient of $E$ of degree $d$ -- this imposes strong conditions on the points $(C,f)\in \overline{\mathcal{M}}_g (\GG (k, E), \beta_d)$: there 
must be one component $C_0$ of $C$ mapping isomorphically onto $X$, while all the other components must be rational, forming trees hanging off the main component $C_0$. The stability condition implies that no such tree can be 
collapsed by $f$. When there are such extra rational components, i.e. when the point is in the \emph{boundary} of $\overline{\mathcal{M}}_g (\GG (k, E), \beta_d)$, the main component $C_0$ gives a section of $\GG (k, E)\rightarrow X$
corresponding to a quotient of $E$ of degree strictly smaller than $d$. Heuristically, 
on one hand such boundary loci have to be rather large because of Bend-and-Break, while on the other hand by the observation above their general points are in the image of glueing (along the marked points) maps of the form
$$\overline{\mathcal{M}}_{g,1} (\GG (k, E), \beta_{d-\delta}) \times_{\GG(k,E)} \overline{\mathcal{M}}_{0,1} (\GG (k, E), \alpha_{\delta}) \longrightarrow \overline{\mathcal{M}}_g (\GG (k, E), \beta_d).$$
Since $d-\delta$ is strictly smaller than $d$, the dimension of the image can be bounded by induction on the degree. In this argument, considering the compactification $\overline{\mathcal{M}}_g (\GG (k, E), \beta_d)$, with ``larger" boundary than $\Quot$, 
is crucial.

I conclude by mentioning a result which is not directly relevant to the sequel, but which is perhaps the main application of Theorem \ref{dimension} outside the realm of moduli spaces. 

\begin{theorem}[\cite{pr}]\label{irreducibility}
For any vector bundle $E$ on $X$, there is an integer $d_Q = d_Q(E,k)$ such that
for all $d\geq d_Q$, ${\rm Quot}_{k,d}(E)$ is irreducible. For any such $d$,
${\rm Quot}_{k,d}(E)$ is generically smooth, has dimension $rd-ke-k(r-k)(g-1)$, 
and its generic point corresponds to a vector bundle quotient.
\end{theorem}

This had been proved before in the case $E = \OO_X^r$ in \cite{bdw}.
It is amusing to note that, reversing the usual order, the result relies on the existence and 
properties of moduli spaces of stable bundles.
In the case when $E$ is semistable, the integer $d_Q$ in the statement above can be made effective (cf. \cite{pr}, end of \S6). There are plenty of examples 
with small $d$ where the result fails in various ways -- cf. e.g. \cite{popa2} \S2 or 
\cite{pr} \S5 for this and other examples of the behavior of Quot schemes on curves.

\section{Verlinde formula and Strange Duality}

The main purpose of this section is to update the comprehensive presentation in \cite{beauville2} with beautiful recent results of Belkale \cite{belkale1}, \cite{belkale2} and Marian-Oprea \cite{mo1}, \cite{mo2} (cf. also \cite{mo3}). 

\subsection{Verlinde formula}

Theorem \ref{picard} says that the only line bundles on the moduli space $SU_X(r, L)$ are the powers of the determinant line bundle $\LL$. The Verlinde formula tells us the dimension $s_{r, d, k}$ of the space of 
global sections $H^0 (SU_X(r, L), \LL^k)$, for all $k \ge 1$, where $d = {\rm deg}(L)$. For an extended discussion of the Verlinde formula and various proofs, see \cite{beauville2}; I will mention two more recent intersection theoretic proofs below. I state a numerical result only in the case when $L = \OO_X$, using the notation $SU_X(r) : = SU_X (r, \OO_X)$ and $s_{r,k} := s_{r, 0, k}$; there are analogous but technically more complicated results for any $L$. In this case, in a form due to Zagier, according to \cite{beauville1} \S5 the formula reads
\begin{equation}\label{verlinde}
s_{r,k} = h^0 (SU_X(r), \LL^k) =  \big(\frac{r}{r+ k}\big)^g \cdot \underset{\underset{|S| = r, |T|= k}{S\cup T = \{1, \ldots, r+k\}}} \sum 
\underset{\underset{t\in T}{s\in S}}{\prod}~ \big{|} 2\cdot {\rm sin}~\pi\frac{s-t}{r +k} \big{|} ^{g-1}.
\end{equation}

\noindent
\begin{example}
Note that $s_{1,k} = 1$, while $s_{r, 1} = r^g$. This last number is the same as the dimension of the 
space of classical theta functions of level $r$, i.e. $h^0 (J(X) , \OO_J (r\Theta))$. We will see in a second that this is not  a coincidence. 
\end{example} 

By Theorem \ref{picard}(4) and Kodaira vanishing (which works on varieties with rational singularities
by Grauert-Riemenschneider vanishing), for any $L$ we have 
$$H^i (SU_X(r, L), \LL^k ) = 0, {\rm ~for~all~} i > 0 {\rm ~and~} k >0,$$
which implies that 
$$h^0 (SU_X(r, L), \LL^k ) = \chi (\LL^k) = \frac{\LL^{(r^2 -1)(g-1)}}{(r^2 -1)(g-1)!}\cdot k^{(r^2 -1)(g-1)} + {\rm l.o.t}.$$
This is a polynomial in $k$, with integral output\footnote{Neither of which seems to be obvious from the formula in the case $L = \OO_X$.} and the knowledge of various intersection numbers on the moduli space appearing in the Riemann-Roch formula naturally leads to the answer. This type of direct calculation proof was first completed by Jeffrey-Kirwan \cite{jk} in the case $(r,d)= 1$ (when the moduli space is smooth), using symplectic methods. Recently they achieved the same in the general case, together with coauthors \cite{jkkw}, using their previous methods in the framework of intersection homology. 

From the algebro-geometric point of view, a surprisingly simple proof was recently provided by Marian-Oprea \cite{mo1} in the smooth case. Hopefully the same type of proof -- with the technical modifications due to the presence of singularities -- will eventually work also in the non-coprime case. 
The point of their proof is to turn the computation into looking at intersection numbers of geometrically defined classes on a Quot scheme parametrizing quotients of a trivial vector bundle. More precisely, let $Q = {\rm Quot}_{r,d} (\OO_X^{r(k+1)})$ be the Quot scheme of coherent quotients of rank $r$ and degree $d$ 
of the trivial bundle of rank $r(k+1)$ on $X$, which can be endowed with a virtual fundamental cycle $[Q]^{vir}$. On $Q \times X$ 
there is a universal bundle $\E$, and the K\"unneth decomposition of the Chern classes of its dual defines the $a_i$ (and $b_i^j$, $f_i$) classes:
$$c_i (\E^\vee) = a_i \otimes 1 + \sum_{j=1}^{2g} b_i^j \otimes \delta_j + f_i, ~i=1,\ldots,r,$$
where $\delta_i$ are generators of the first cohomology of $X$. 
It turns out that only some very special intersection numbers on this Quot scheme are needed in order to obtain the Verlinde formula. In fact the key formula proved in \cite{mo1} is
\begin{equation}\label{intersection}
h^0 (SU_X(r, L), \LL^k ) = \frac{1}{(k+1)^g} \cdot \int_{[Q]^{vir}} a_r^{k( d - r(g-1)) + d}.
\end{equation}
Now, at least when $d$ is large enough, this intersection number has an enumerative meaning, namely it counts the number of 
degree $d$ maps to the Grassmannian $\GG (r, r(k+1))$ with incidences at $k( d - r(g-1)) + d$ points with sub-Grassmannians 
$\GG (r, r(k+1)-1)$ in general position.

What about $U_X(r,d)$? For instance, note that we have a natural morphism obtained by taking tensor product with line bundles:
$$\tau : SU_X(r) \times {\rm Pic}^{g-1}(X) \overset{\otimes}{\longrightarrow} U_X(r, r(g-1)).$$
One can check by a slight variation of Exercise \ref{pullback} that if we denote the canonical polarization 
$$\Theta = \Theta_{can} : = \{ E ~|~ h^0 E \neq 0 \} \subset U_X(r, r(g-1)),$$
we have that 
$$\tau^* \OO_U (\Theta) \cong \LL \boxtimes \OO_{{\rm Pic}}(r \Theta).$$
The K\"unneth formula gives then for any $k \ge 1$ an isomorphism
$$H^0 (\tau^*\OO_U (k\Theta)) \cong H^0 (SU_X(r), \LL^k ) \otimes 
H^0 ({\rm Pic}^{g-1}(X), \OO(r \Theta)).$$
Putting these facts together, we obtain that 
$$h^0 (U_X(r, r(g-1), \OO_U (k\Theta)) = s_{r,k} \cdot \frac{k^g}{r^g}.$$

\begin{corollary}\label{U_dimension}
$$h^0 (U_X(k, k(g-1), \OO_U (r\Theta)) = s_{r,k} = h^0 (SU_X(r), \LL^k).$$
\end{corollary}
\begin{proof}
Stare at the Verlinde formula until you see that the expression $s_{r,k}\cdot r^{-g}$ is 
symmetric in $r$ and $k$.
\end{proof}

\begin{corollary}\label{principal}
$$h^0 (U_X(k, k(g-1), \OO_U (\Theta)) = 1, {\rm~for~all~} k \ge 1.$$
\end{corollary}

More generally, in \cite{dt} and \cite{beauville1} \S8 it is explained that a simple covering argument reduces the calculation of
the dimension of spaces of generalized theta functions on $U_X(r,d)$ to the Verlinde formula as well. 
Here, as before, $\Theta_F$ is any \emph{basic} theta divisor on $U_X(r,d)$ as in Theorem \ref{picard}(3), $h = (r,d)$, and $s_{r,d,k} = h^0 (SU_X(r,L), \LL^k)$ for any $L\in {\rm Pic}^d(X)$.

\begin{proposition}
$$h^0 (U_X(r, d), \OO_U (k \Theta_F)) = s_{r,d,k}\cdot\frac{k^g}{h^g}.$$
\end{proposition}

\subsection{Strange Duality}
Furthermore, recalling the notation $r_0 = r/h$ and $d_0 = d/h$, one defines a tensor product map: 
$$\tau : SU_X(r,L) \times U_X (kr_0, k(r_0(g-1)-d_0)) \overset{\otimes}{\longrightarrow} U_X(krr_0, krr_0(g-1)).$$
We pull back the unique section given by Corollary \ref{principal} to the product space. Hence up to constants we have a naturally defined section in $H^0 (\tau ^* \OO_U (\Theta))$. Consider a stable bundle $G$ of rank $r_0$ and degree $d_0$
which has a theta divisor. It can be checked easily that this is a basic theta divisor on $U_X (kr_0, k(r_0(g-1)-d_0))$
If we impose the (harmless) assumption that $L \cong ({\rm det}~G)^{\otimes h}$, we have a simple formula for the pull-back line bundle: 
$$\tau^*\OO_U (\Theta)\cong \LL^k \boxtimes \OO_U (h\Theta_G).$$
Thus our section lives naturally in 
$$H^0(SU_X(r,L),\LL^k)\otimes H^0(U_X(kr_0,k(r_0(g-1)-d_0)), \OO_U(h\Theta_G)),$$ 
and so it induces a non-zero homomorphism 
$$SD: H^0(SU_X(r,L), \LL^k)^\vee \longrightarrow H^0(U_X(kr_0,k(r_0(g-1)-d_0), \OO_U(h\Theta_G)).$$
We saw that the two spaces have the same dimension. In fact much more is true, namely what used 
to be called the \emph{Strange Duality Conjecture}.

\begin{theorem}\label{sd}
$SD$ is an isomorphism. 
\end{theorem}

This is certainly one of the most exciting recent theorems to be proved in the context of generalized theta functions. Belkale \cite{belkale1} first proved it for general curves. Subsequently, building on some of his ideas, and introducing several new ones as well, 
Marian and Oprea \cite{mo2} proved it for all curves.  Finally, in \cite{belkale2}, Belkale complemented the picture by showing the flatness of the global Strange Duality map over the moduli space of curves. A more detailed discussion is given below.
(Previous results include a proof for $k =1$ by Beauville-Narasimhan-Ramanan \cite{bnr}, and a proof for $r = 2$ by Abe \cite{abe}.)

\begin{remark}\label{interpretation}
The Strange Duality has a beautiful geometric interpretation (which is also the key to its proof): \emph{for any $k \ge 1$, the linear series $|\LL^k|$ on $SU_X(r,L)$ is \emph{spanned} by the generalized theta divisors $\Theta_F$, with $F \in U_X(kr_0, k(r_0(g-1)-d_0))$}. Here is a sketch of this, following \cite{beauville1} \S8. The duality is equivalent to having a commutative diagram
$$\xymatrix{
SU_X(r,L)  \ar[r]^{\phi_{\LL^k}} \ar[dr]^{\theta_k} &  {\bf P}((H^0 \LL^k)^\vee)  \ar[d]_{\cong}^{SD} \\
& {\bf P}(H^0 \OO_U (h\Theta_G)) } $$ 
where $\phi_{\LL^k}$ and $\theta_k$ are in general only \emph{rational} maps, the first one induced by the line bundle $\LL^k$, and the second taking the form $E \to \Theta_E \in |h\Theta_G|$. Note that the divisor of the section inducing the $SD$ map is 
$$\mathcal{D} = \{ (E, F) ~|~ h^0 (E \otimes F) \neq 0 \} \subset SU_X(r,L) \times U_X(kr_0,k(r_0(g-1)-d_0).$$
If $H_F \in {\bf P}(H^0 \OO_U (h\Theta_G))$ is the hyperplane of divisors passing through $F$, then 
$\theta_k^* H_F = \Theta_F$, so via this duality the generalized theta divisors $\Theta_F$ correspond 
to the hyperplanes $H_F$, which obviously span the linear series  ${\bf P}(H^0 \OO_U (h\Theta_G))$.
\end{remark}

I will give a quick overview of the methods involved in the proof of the Strange Duality, and for technical simplicity I stick to the case $d=0$. Belkale's approach \cite{belkale1} is to use the geometric interpretation 
of the Strange Duality in the Remark above. Namely, his aim is to construct $m = s_{r,k}$ vector bundles $E_1, \ldots, E_m$ with trivial 
determinant, and also $m$ vector bundles $F_1, \ldots, F_m$ of rank $k$ and degree $k (g-1)$, such that 
$$H^0 (E_i \otimes F_j) = 0 \iff i = j.$$
This implies that there exist $m$ linearly independent generalized theta divisors $\Theta_{F_1}, \ldots, \Theta_{F_m}$ in the linear
series $|\LL^k|$ of dimension $m$, hence the result. Belkale achieves this for the general curve $X$ of genus $g$ by degeneration to a rational curve with $g$ nodes. On such a curve, the problem is translated into an explicit Schubert calculus problem for a product of Grasmannians. He then uses 
the connection between the quantum cohomology of Grassmanians and representation theory, in particular the Verlinde formula, known already from work of Witten \cite{witten}. The subtlety in the argument comes from choosing a certain part of the intersection numbers involved, which corresponds
to fixing the determinant of the $E_i$. 

To obtain the result on any curve, Marian and Oprea \cite{mo2} follow Belkale's strategy, noting in addition that it is not the Verlinde numbers $s_{r,k}$, but rather the numbers $\frac{(r+k)^g}{r^g} \cdot s_{r,k}$ that have enumerative meaning. More precisely, in this case formula $(\ref{intersection})$ can be rewritten as
$$h^0 (SU_X(r), \LL^k ) = \frac{r^g}{(r+k)^g}\cdot \int_{Q^\prime} a_k^s,$$
where this time $Q^\prime$ stands for the Quot scheme of quotients of a trivial bundle of rank $r+k$.
Here $s:= (r+k) \frac{d}{k} - r(g-1)$, $d$ is taken to be divisible by $k$, and sufficiently large so that in particular $Q^\prime$ is irreducible and of expected dimension (cf. \cite{bdw}). As noted before, for such $d$ the intersection numbers $\int_{Q^\prime} a_k^s$ count maps from rational curves to Grassmannians, with certain incidence conditions. After some standard manipulation, this implies that they count vector bundles $E_i$ as above, without the determinant restriction. 
On the other hand, a slight variation on Corollary \ref{principal} shows that they represent the dimension of certain 
spaces of sections on $U_X(r,0)$. Marian and Oprea formulate a new duality statement for these moduli spaces, which can be seen as an extension of the Wirtinger duality for $2\Theta$ functions on Jacobians. The count and interpretation above show precisely the following:

\begin{theorem}(\cite{mo2})
Let $L$ and $M$ be two line bundles of degree $g-1$ on $X$. There exists a rank-level duality isomorphism
$$D: H^0 (U_X(r,0),\OO_U(k \Theta_M)\otimes {\rm det}^*\OO_J(\Theta_L))^\vee \longrightarrow$$  
$$\longrightarrow 
H^0 (U_X(k,0), \OO_U(r\Theta_M) \otimes {\rm det}^* \OO_J(\Theta_{\omega_X \otimes L^{-1}})).$$
\end{theorem}

A restriction argument gives from this without much effort the usual Strange Duality statement. The argument for arbitrary degree is very similar. Finally, note that recently Belkale \cite{belkale2} has shown that the \emph{global} Strange Duality map over the moduli space of curves is projectively flat (cf. also \cite{beauville2} for background on 
this circle of ideas), which beautifully completes the picture. This says in particular that knowing the Strange Duality for any one curve is equivalent to knowing it for an arbitrary curve.

\section{Base points}

\subsection{Abstract criteria}

The geometric interpretation of the Strange Duality stated in Remark \ref{interpretation} provides a simple 
criterion for detecting base points for any linear system $|\LL^k|$ on an arbitrary $SU_X(r, L)$  
with $L \in {\rm Pic}^d(X)$.

\begin{lemma}\label{base_point_criterion}
A vector bundle (class) $E$ in $SU_X(r, L)$ is a base point for the linear system $|\LL^k|$
if and only if 
$$H^0 (E\otimes F) \neq 0, ~{\rm for~ all}~ F\in U_X(kr_0, k(r_0(g-1) -d_0)).$$
\end{lemma}  

Recently, another abstract description of the base points of these linear series has emerged. I will briefly describe the general construction. 
It is well known that there exists an \'etale cover $M$ of $U_X(kr_0, k(r_0 (g-1) -d_0))$ such that there is a ``universal bundle" $\E$ on $X\times M$ (cf. \cite{nr2} Proposition 2.4). This allows one to consider a Fourier-Mukai-type correspondence ${\bf R} \Psi_{\E}: {\bf D}(M) \rightarrow {\bf D}(X)$ taking an object $A$ to ${\bf R}{p_X}_* (p_M^* A \otimes \E)$.\footnote{For the statements we are interested in, this is a good as thinking that $M$ is $U_X(kr_0, k(r_0 (g-1) -d_0))$ itself, with the technical problem that for $k \ge 2$, and sometimes for $k=1$, this moduli space is not fine.} 
We fix any generalized theta divisor (or indeed any ample line bundle) on $U_X(kr_0, k(r_0 (g-1) -d_0))$, and denote abusively by $\Theta$ its pullback to $M$. We consider for $m \gg 0$ the Fourier-Mukai transform
$$E_m^k : = {\bf R} \Psi_{\E} \OO_M (-m\Theta) [\dim M] = R^{{\rm dim}~M} \Psi_{\E} \OO_M(-m\Theta),$$
which is a vector bundle on $X$.\footnote{This, and the second equality above, are due to the fact that $m\Theta$ is 
taken sufficiently positive. We abusively omit in the notation the choice of ample line bundle made on $M$.}
The equivalence between (1) and (3) in the statement below was first proved by Hein \cite{hein2} (noting that objects 
essentially equivalent to the vector bundles $E_m^k$ can be defined alternatively in terms of spectral covers of $X$).
It follows also, together with the equivalence with (2), as a special case of the theory of $GV$-sheaves in \cite{pp2}.

\begin{theorem}(\cite{hein2} and \cite{pp2})\label{hein}
Let $X$ be a smooth projective curve of genus $g \ge 2$ and $E$ a vector bundle in
$SU_X(r, L)$, with $L \in {\rm Pic}^{d}(X)$. Then the following are equivalent:
\begin{enumerate}
\item $E$ is not a base point for the linear series $|\LL^k|$.
\item $H^1 (E \otimes E_m^k) = 0$ for all $m \gg 0$.
\item $H^0 (E \otimes {E_m^k}^\vee) = 0$ for all $m \gg 0$.
\end{enumerate}
\end{theorem}

In the most basic, and by far most studied case, namely when 
$d= r(g-1)$ (equivalent by appropriate tensorization to the case $d=0$) and $k = 1$, the Fourier-Mukai 
correspondence above is between $X$ and ${\rm Pic}^0(X)$, and $\Theta$ is simply a theta divisor on 
the Picard variety of $X$. Here the bundles $E_m^1$ are well understood: they are precisely the Raynaud bundles described in \S2.4, which as we saw are stable. Moreover, Hein \cite{hein2} provides an effective bound in (3), which gives hope for a more or less concrete description of the entire base locus. 
It would be nice to similarly have a better understanding 
of \emph{all} the vector bundles $E_m^k$ considered in Theorem \ref{hein}. For example one would expect:

\begin{conjecture}
All the vector bundles $E_m^k$ are stable.
\end{conjecture}

\subsection{Concrete constructions}

Lemma \ref{base_point_criterion} was until recently conjectural in essentially every other 
case except for $L = \OO_X$ and $k = 1$, so most of the existing constructions of base points are 
in this case.  Putting together what is known up to now, in this case the conclusion is roughly the following:

\noindent
(1) There exist basic examples of base points, and even basic families of large dimension in special ranks, obtained both from the Raynaud examples \cite{raynaud}, and from the examples in \cite{popa1}. Still basic 
examples, following the ideas of the previous two, are examples in \cite{popa3} (see also Remark
\ref{verlinde_other} (2)), respectively \cite{schneider1} and \cite{gt}. Other examples can be built by applying standard operations (extensions and elementary transformations) to the basic examples, as in \cite{arcara} or \cite{schneider2},
and similarly, but more subtly, in \cite{pauly}, based also on some ideas in \cite{hitching}.

\begin{example} 
Let me explain the most basic constructions.
To begin with, let's restate Corollary \ref{base_point_criterion} in this case.

\begin{corollary}
$E$ is a base point for $|\LL|$ if and only if $H^0 (E \otimes L ) \neq 0$ for all $L 
\in {\rm Pic}^{g-1}(X)$. 
\end{corollary} 

By twisting with an effective line bundle in order to obtain the slope $\mu = g-1$, we see that it 
is enough to construct semistable vector bundles $E$ satisfying 
$$ (*) ~ \mu (E) \in \ZZ, ~ \mu (E) \le g-1, h^0 ( E \otimes \xi ) \neq 0 ~~ {\rm for~all~}\xi \in {\rm Pic }^0 (X).$$

\noindent 
{\bf Lazarsfeld's bundles II \cite{popa1}.}
Let's go back to the Lazarsfeld bundles considered in \S2.3. 
Let $L$ be a line bundle on $X$ of degree $d\geq 2g+1$ and consider the evaluation sequence 
$$0 \longrightarrow M_{L} \longrightarrow H^{0}(L)\otimes \mathcal{O}_{X}
\longrightarrow L \longrightarrow 0.$$
We saw that $Q_{L}= M_{L}^{\vee}$ is a stable bundle, and $\mu(Q_L) = \frac{d}{d-g}< 2$.
Recall that we showed that if $x_1, \ldots, x_{d-g-1}$ are general points on $X$, then we have an 
inclusion
$$ \bigoplus_{i=1}^{d-g-1}\mathcal{O}_{X}(x_{i})
\hookrightarrow Q_{L}.$$
Passing to exterior powers, we obtain for any $p$ an inclusion
$$ \underset{1\le j_1 < \ldots < j_{p}\le d-g-1}{\bigoplus}\OO_{X}(x_{j_1} + \ldots  + x_{j_{p}})
\hookrightarrow \bigwedge^p Q_L.$$
This has the consequence that
$$h^0 (\bigwedge^p Q_L \otimes \OO_X(y_1 + \ldots +y_p - x_1 -\ldots - x_p)) \neq 0$$
for any general points $y_1, \ldots, y_p$ and $x_1, \ldots, x_p$ on $X$. Such differences of points cover what is called the \emph{difference variety} $X_p - X_p \subset {\rm Pic}^0 (X)$. This has 
dimension $2p$, hence in particular we have ${\rm Pic}^0 (X)= X_{\gamma} - 
X_{\gamma}$, where $\gamma := [\frac{g+1}{2}]$.  By semicontinuity we finally deduce that 
$$h^0 (\bigwedge^{\gamma} Q_L \otimes \xi) \neq 0, {\rm~for~all~} \xi \in {\rm Pic}^0(X).$$
Since at least in characteristic zero $\bigwedge^{\gamma} Q_L$ is semistable,  in order to satisfy 
condition $(*)$ we only need to make a choice such that $\mu(\bigwedge^{\gamma} Q_L) = \gamma \cdot \mu(Q_L) = \frac{\gamma d}{d-g} \in \ZZ$. This is certainly possible: the simplest such choice is to take $d = g (\gamma + 1)$, which gives $\mu(\bigwedge^{\gamma} Q_L) = \gamma + 1$.
In conclusion, starting with a sufficiently high rank, $|\LL|$ will have base points.\footnote{Note that starting with a base point in rank $r$, one can always obtain a base point in any rank $s > r$ by adding  
for example the trivial bundle of rank $s-r$.} It is explained in \cite{popa1} that a more careful construction provides higher dimensional families of such examples (with fixed determinant).

\noindent
{\bf Raynaud's bundles II \cite{raynaud}.} 
Raynaud's vector bundles discussed in \S2.4 also provide such examples. This was the first known construction of base points for $|\LL|$, as noted by Beauville \cite{beauville2}.
Let's see how it works: recall that such a bundle $E$ is 
the restriction to an Abel-Jacobi embedding $X \subset J(X)$ of the Fourier-Mukai transform $F = \widehat{\OO_{J(X)} (-m\Theta)}$.

Consider $U \subset J(X)$ an open subset where $\OO_{J(X)} (m\Theta)$ is trivial. By Mukai's inversion theorem \cite{mukai2} Theorem 2.2, $\widehat F \cong (-1_{J(X)})^* \OO_{J(X)} (- m\Theta)$, 
so it is trivial on $U$ as well. By definition, 
$\widehat F \cong {p_2}_*(p_1^* F \otimes \mathcal{P})$. Thus our choice means that there exists 
$0 \neq s \in \Gamma (p_1^* F \otimes \mathcal{P} _{|p_2^{-1}(U)}).$
Hence there exists $x_0 \in J(X)$ such that $s_{| \{x_0\} \times U} \neq 0$. But we can choose our Abel-Jacobi embedding of $X$ to pass through $x_0$. This implies that $s_{|X\times U} \in \Gamma (p_1^*F \otimes \mathcal{P}_{| X\times U})$ is non-zero, which reinterpreted means that $H^0 (E \otimes \xi) \neq 0$ for $\xi \in {\rm Pic}^0(X)$ general. By semicontinuity this has to hold then for all $\xi \in {\rm Pic}^0(X)$. 
As we saw in  \S2.4, the slope of $E$ is $\frac{g}{m}$. If $m|g$, then $E$ satisfies condition $(*)$.
It turns out that these examples are the first level of a more general class of bundles with the same 
property, of rank $m^g$ and slope $\frac{kg}{m}$, for any $k< m$ coprime with $m$; see    
Remark \ref{verlinde_other} (2).
\end{example}

\noindent
(2) For each genus $g$, there exist base points for the linear series $|\LL|$ on $SU_X(r)$
for every rank $r \ge r_0 = g+2$, while for hyperelliptic curves they exist for every rank $r \ge r_0 = 4$ 
(see \cite{pauly}; cf. also \cite{arcara} for the weaker bound $r \ge r_0 = \frac{1}{2} (g+1)(g+2)$). This 
suggests that the lower bound $r_0$ might depend on the curve $X$.

\noindent
(3) If $B_r : = {\rm Bs}~|\LL| \subset SU_X(r)$, then 
$${\rm dim}~ B_r  + (s-r)^2 (g-1) + 1 \le {\rm dim}~ B_s, ~{\rm ~for~all~} s\ge r.$$
This follows roughly speaking by adding all semistable bundles of rank $s-r$ to the base points (cf. \cite{arcara}).

My favorite problem remains giving examples of basic base points, as in (1). It is probably no accident that essentially all nontrivial classes of examples of vector bundles on curves we know of have to do with this problem! Also, given the successive improvements, it is perhaps time to ask:

\begin{question}
What is a good estimate for the dimension of the base locus of $|\LL|$? Moreover, can one give a reasonably concrete description of the entire base locus?
\end{question}

A full description of the base locus is known in genus $2$ and rank $4$: it is a reduced finite set of cardinality $16$, consisting of Raynaud-type bundles (cf. \cite{hitching} and \cite{pauly}).

\subsection{Other determinants or levels}

A base point for $|\LL^k|$, $k \ge 2$, is of course also a base point for $|\LL|$, so the most comprehensive construction problem is that of base points of $|\LL|$. In \cite{popa1} it is explained however that there is a concrete way of constructing base points for higher $k$ (and also for determinants $L\neq \OO_X$). This is an indirect construction, based on combining Raynaud's examples and the existence of subbundles of large enough degree of a fixed 
vector bundle. 

\begin{lemma}\label{subbundle}
Let $E$ be a semistable bundle on $X$, such that $h^0 (E \otimes \xi ) \ne 0$ for all $\xi \in {\rm Pic}^0(X)$. If $k \le \frac{g}{\mu(E) + 1}$ is an integer, then 
$H^0 (E \otimes F) \ne 0$ for all vector bundles $F$ with $\mu(F) = g-1-\mu(E)$.
\end{lemma}
\begin{proof}
If $F$ has rank $k$ and degree $e$, then by the theorem of Lange and Mukai-Sakai (cf. Remark \ref{subbundles}(3)) we know that there exists a sub-line bundle $M\hookrightarrow F$ of degree 
$${\rm deg}(M)\geq \frac{e-(k-1)g}{k}.$$
The inequality in the hypothesis implies precisely that ${\rm deg}(M) \ge 0$, which, given the other 
condition in the hypothesis, implies that $H^0 (E \otimes M) \neq 0$. This gives immediately that 
$H^0 (E \otimes F) \neq 0$.
\end{proof}

\noindent
Combining this with Lemma \ref{base_point_criterion}, one obtains

\begin{corollary}
Let $L \in  {\rm Pic}^d(X)$ for some $d$. If $E$ is a vector bundle in $SU_X(r,L)$ satisfying Raynaud's condition:
$$H^0 (E \otimes \xi) \neq 0, {\rm ~for~all~} \xi \in {\rm Pic}^0 (X),$$
then $E \in {\rm Bs}|\LL^k| $ for all integers $k\le \frac{g}{d_0 + r_0}$.
\end{corollary}

\noindent
One concrete consequence that deserves special notice is 

\begin{corollary}
Raynaud's examples with integral slope $\mu = g/n$ produce, after twisting with appropriate line bundles, base points in $|\LL^k|$ on $SU_X(r)$ if $k\le \frac{gn}{g+n}$, where $r = n^g$. 
\end{corollary}

Thus $k$ can be made arbitrarily large if $g$ and $n$ are chosen appropriately, which rules out the possibility of a fixed bound on all $SU_X(r)$. The reader can draw similar conclusions in other degrees by playing with the known examples satisfying Raynaud's condition, listed in the previous subsection.

\begin{problem}
Find a direct way of constructing vector bundles $E$ satisfying the condition $H^0 (E\otimes F) \neq 0$ for all semistable bundles $F$ of slope $g- 1 -\mu(E)$ and some rank $k \ge 2$. 
\end{problem}

One construction proposed in \cite{popa1}  (Remark in \S2) is to consider vector bundles of the form 
$\bigwedge^p Q_E$, analogous to the Lazarsfeld bundles above, where this time $E$ is a higher rank vector bundle as opposed to a line bundle. This was carried out recently in 
\cite{cgt} for arbitrary determinants (cf. also \cite{gt}), providing positive dimensional families of base points in various linear series on all moduli spaces.

\section{Effective results on generalized theta linear series}

\subsection{The theta map}

Not much is known at the moment about the rational map 
$$\theta: SU_X(r) \dashrightarrow |r\Theta|, ~~E \to \Theta_E,$$ 
which by Strange Duality can be identified with the map induced by the linear series $|\LL|$. This is the most natural map defined on the moduli space, but unfortunately we have already seen in the previous section that usually it is not a morphism. A very nice survey of what is currently known on the theta map is \cite{beauville5}. The reader is referred to it, while here I will only briefly reproduce the main results. 

\begin{theorem}\label{morphism}
(1) If $r=2$, then  $\theta$ is a morphism (\cite{raynaud}; cf. also \cite{hein1} and \cite{popa2}).

\noindent
(2) Moreover, if $r=2$ and $X$ is non-hyperelliptic, then $\theta$ is an embedding (\cite{bv1} and \cite{vgi}), while if $X$ is 
hyperelliptic, it is $2:1$ onto its image (cf. \cite{beauville1}).

\noindent
(3) If $r=3$, then if $X$ is general $\theta$ is a morphism (\cite{raynaud}). The same is true for all $X$ if in addition $g = 2$ (\cite{raynaud}) or $g = 3$ (\cite{beauville4}).

\noindent
(4) If $g = 2$, $\theta$ is generically finite (\cite{beauville4}).
\end{theorem}

In a few instances the image of the theta map has been determined explicitly. This is summarized below. It remains a fascinating problem to determine explicitly the moduli space of vector bundles as an embedded projective variety in other situations. 

\begin{theorem}
(1) If $g=2$, then $SU_X(2) \cong \PP^3$ (\cite{nr1}).

\noindent
(2) If $g=3$ and $X$ is non-hyperelliptic, then $SU_X(2)$ is isomorphic to the Coble 
quartic hypersurface in $\PP^7$ (\cite{nr3}).

\noindent
(3) If $g=2$, then $SU_X(3)$ is a double covering of its image in $\PP^8$, branched along  a sextic which is the dual 
of the Coble cubic hypersurface (\cite{ortega} and \cite{nguyen}).
\end{theorem}

Even less studied is the theta-map for other $SU_X(r,d)$, besides the case of rank $2$. Here are some 
results that I am aware of. 

\begin{theorem}
(1)  $|\LL|$ is very ample on $SU_X(2,1)$ (\cite{bv2}).

\noindent
(2) If $X$ is hyperelliptic of genus $g$, then $SU_X(2,1)$ is isomorphic to the subvariety of the 
Grassmannian ${\bf G}(g-2, 2g+1)$ parametrizing the $(g-2)$-planes contained in the intersection 
of two general quadrics in $\PP^{2g+1}$ (\cite{dr}).

\noindent
(3) For arbitrary $r$, the rational theta map on $SU_X(r,1)$ is generically injective (\cite{bv3}).\footnote{In \emph{loc. cit.} it is also stated that it is a morphism, but unfortunately this seems to have a gap.}
\end{theorem}

One last interesting thing to note is Beauville's observation \cite{beauville3} that there exist (even stable) examples of vector bundles $E$, at least of integral slope, which are not characterized by their theta divisor $\Theta_E$. In other words, for every $r\ge 3$, there exist genera for which the theta map 
on $SU_X(r)$ is not injective.

\subsection{Basepoint-freeness} 

Since the linear series $|\LL|$ is most of the time not basepoint-free on $SU_X(r,d)$, the 
next natural thing is to look for effective bounds for the basepoint-freeness of powers $|\LL^m|$. The remarks at the end of \S6 suggest that one should not expect a bound independent of the 
rank $r$. On the other hand, Faltings has shown in \cite{faltings} that one could eventually generate this linear series only by generalized theta divisors, without providing an effective bound. 

Effective versions of Faltings' results are now available in \cite{popa2} and \cite{pr}, improving on earlier results of Le Potier \cite{lepotier1} and Hein \cite{hein1} which depended on the genus of the curve. In order to state the general result, I will need some notation. Let $d$ be an arbitrary integer and let as usual $h= (r,d)$. For the statement it is convenient to introduce another invariant of the moduli space. If $E\in SU_{X}(r,d)$ and $1\leq k\leq r-1$, as in \cite{lange} we define
$s_{k}(E):=kd-r e_{k}$, where $e_{k}$ is the maximum degree of a 
subbundle of $E$ of rank $k$. Note that if $E$
is stable, one has $s_{k}(E)\geq h$ and we can further define 
$$s_{k}=s_{k}(r,d):=\underset{E~ {\rm stable}}{{\rm min}}s_{k}(E)
{\rm ~and~}  s=s(r,d):=\underset{1\leq k\leq r-1}{{\rm min}}s_{k}.$$ 
Clearly $s\geq h$ and it is also immediate that $s(r,d)=s(r,-d)$.

\begin{theorem}\label{bpf}
(1) (\cite{popa2}) The linear series $|\LL^m|$ on $SU_{X}(r,d)$ is basepoint-free
for $$m\geq {\rm  max}\{[\frac{(r+1)^{2}}{4r}] h,[\frac{r^{2}}{4}]\frac{h}{s}\}.$$ 

\noindent
(2) (\cite{pr}) Slightly better, on $SU_X(r)$ the linear series $|\mathcal{L}^{m}|$ is basepoint-free  for 
$m\geq [\frac{r^{2}}{4}]$.
\end{theorem}

The general pattern is that for arbitrary degree one has bounds which are at least as good as those for degree $0$. For example, if $d = \pm 1$, then $m \ge r-1$ suffices.

\begin{remark}
What should we have definitely expected at worst? Fujita's Conjecture predicts basepoint-freeness for $\omega_{SU_X} \otimes \LL^p$ with $p \ge {\rm dim} ~SU_X + 1$ (and one more for very ampleness). Given the formula for $\omega_{SU_X}$ in Theorem \ref{picard}, this would predict the bound $k \ge (r^2 -1) (g-1) + 1 -2h$ for basepoint-freeness, with $h = (r, d)$, which is also quadratic in $r$ but depends on the genus.
\end{remark}

\noindent
\emph{Sketch of proof of Theorem \ref{bpf} (2) (cf. \cite{popa2} \S4 for details).} (Part (1) goes along the same lines, but is slightly more technical.) The crucial point is to relate this to Quot schemes.

\begin{proposition}
If $E \in SU_X(r)$ and ${\rm dim}~ {\rm Quot}_{k,d}(E) \le m\cdot d$ for all $1 \le k \le r-1$, then 
$E \not \in {\rm Bs}~|\LL^k|$.
\end{proposition} 
\begin{proof}(\emph{Sketch.})
It suffices to show that there exists $F \in U_X(m, m(g-1))$ such that $H^0 (E \otimes F) = 0$. Let's ``count" the moduli of $F$ for which this is not the case. In other words, assume that $H^0 (E \otimes F) \neq 0$, i.e. that there exists a non-zero homomorphism $E^\vee \overset{\phi}{\rightarrow} F$. This induces a diagram of the form
$$\xymatrix{
E^\vee \ar[r] \ar[dr]_{\phi} &  G \ar[r] \ar@{^{(}->}[d] & 0 \\
 & F }$$
where the vector bundle $G$ is the image of $\phi$.
The vertical inclusion determines an exact sequence
\begin{equation}
0\longrightarrow G\longrightarrow F\longrightarrow H^{\prime}\longrightarrow 0,
\end{equation}
where $H^{\prime}$ might have a torsion part $\tau_{a}$ of length $a$. This induces a diagram as follows, where $H$ is locally free of rank $m-k$ and degree $m(g-1) - d- a$.

$$\xymatrix{
& & 0 \ar[d] & 0 \ar[d] \\
0 \ar[r] & G \ar[r] \ar[d]_{\cong} & G^{\prime} \ar[r] \ar[d] & \tau_{a} \ar[r]
\ar[d] & 0 \\
0 \ar[r] & G \ar[r] & F \ar[r] \ar[d] & H^{\prime} \ar[r] \ar[d] & 0 \\
& & H \ar[d] \ar[r]^{\cong} & H \ar[d] \\
& & 0 & 0 }$$ 
Clearly then, for fixed degree, the moduli of such $G$'s and $G^{\prime}$'s satisfy
$${\rm dim}~\{G^{\prime}\} \le {\rm dim}~\{G\} + ka.$$
Now note that if $d$ is the degree of 
$G$, then $d$ runs only over a finite set of integers, since $E^\vee$ and $F$ are
semistable. Thus we can do the calculation for any fixed $d$. 
Now the dimension of the moduli of such $F$'s is at most the dimension of the space of vertical extensions in the diagram above, so we can try to estimate that. We have 

\begin{equation}\label{first}
{\rm dim}~\{F\} \le {\rm dim}~\{G\}_d  + ka + {\rm dim}~\{H\}_{m(g-1) -d -a} 
+ {\rm dim} ~{\rm Ext}^1 (H\oplus \tau, G) -1.
\end{equation}

On the other hand, the dimension of the moduli of $G$'s in the diagram can be bounded 
above by the dimension of the space of quotients of $E$ of the type of $G$, which is in turn 
bounded by hypothesis:

\begin{equation}\label{second}
{\rm dim}~\{G\}_d \le {\rm dim} ~{\rm Quot}_{k,d}(E) \le m\cdot d.
\end{equation}

Since $F\in U_X( m, m(g-1))$ form a bounded family, for a fixed degree the $H$'s must form 
a bounded family too (this is in fact a subset of a \emph{relative} Quot scheme). The well-known Lemma \ref{upper_bound} below implies then the estimate

\begin{equation}\label{third}
{\rm dim}~\{H\}_{m(g-1) -d -a} \le ({\rm rk}~H)^2 (g-1) + 1 = (m-k)^2 (g-1) +1.
\end{equation}

\noindent
For the space of extensions, using Serre duality we have 

\begin{equation}\label{fourth}
{\rm dim}~{\rm Ext}^1 (H\oplus \tau_{a}, G) = h^1 (H^\vee \otimes G) + ka.
\end{equation}
Finally, since $F$ is stable, by Exercise \ref{stable_sequence} we have $h^0 (H^\vee \otimes G) = 0$, so we can use Riemann-Roch to compute

\begin{equation}\label{fifth}
h^1 (H^\vee \otimes G) = k(m-k) (g-1 - \frac{d}{k} + \frac {m(g-1) - d- a}{m-k}).
\end{equation}

\noindent
Putting together (\ref{first}), (\ref{second}), (\ref{third}), (\ref{fourth}) and (\ref{fifth}), after a slightly tedious calculation one gets 
$${\rm dim}~\{F\} \le m^2 (g-1) < {\rm dim}~U_X(m, m(g-1)),$$
so for the general $F$ we must have $H^0 (E \otimes F) = 0$.
\end{proof}

\begin{lemma}[cf. \cite{nr2}, Proposition 2.6]\label{upper_bound}
Let $S$ be a scheme of finite type over $k$, and $\E$ on $X\times S$ a vector bundle such that all $\E_s := \E _{|X\times \{s\}}$ have rank $r$ and degree $d$. Then there exists a variety $T$ (irreducible, 
nonsingular), and a vector bundle $\F$ on $X\times T$, such that all vector bundles $\E_s$ and all 
stable bundles of rank $r$ and degree $d$ are among the $\F_t$'s, $t\in T$.
\end{lemma}

Going back to the proof of Theorem \ref{bpf},  now the key point is to use the bound for the dimension of Quot schemes provided by Theorem \ref{dimension}. In view of this and the discussion above, what's left to note is that if $E\in SU_X(r)$ is stable, then 
$$k(r-k) + (d - d_k) r \le d \cdot [\frac{r^2}{4}]$$ 
for all $d$ and all $k$.

\subsection{Very ampleness} 

Aside from the beautiful results of Brivio-Verra and van Geemen-Izadi mentioned above, showing that in rank $2$ the theta map is usually an embedding, in general one has to go to even higher powers than in the previous subsection in order to achieve this. Faltings' basepoint-freeness method was extended by Esteves \cite{esteves} to point separation and separation of tangent vectors at smooth points, the starting point being the criterion in Corollary \ref{effective}(3). This, together with the method in \cite{popa2} described in \S7.2, is used in \cite{ep} to obtain effective very ampleness bounds. There is however an obstruction to obtaining a complete result for arbitrary rank and degree: at the moment we don't know how to separate tangent vectors at singular points of the moduli space using generalized theta divisors. Thus very ampleness is obtained only in the case $(r,d) = 1$, while in general we only have an injective map which separates tangent vectors at smooth points.

\begin{theorem}(\cite{ep})\label{very_ampleness}
For each $m\geq r^2+r$ the linear series $|\LL^m|$ on $SU_X(r,d)$ separates points, and separates tangent vectors on the smooth locus $SU^s_X(r,d)$. In particular, if $r$ and $d$ are coprime, then $|\LL^m|$ is very ample.
\end{theorem}

Note that here, unlike in Theorem \ref{bpf}, the result does not depend on the degree. 
I suspect that the bound for very ampleness should hold in general.\footnote{Note however that a sharper bound is proposed in Conjecture \ref{bounds}.}

\begin{conjecture}
The very ampleness bound $r^2 + r$ holds for any rank and degree.
\end{conjecture} 

In order to attack this problem, one likely needs a good understanding of the tangent cone, and space, to the moduli space of vector bundles at a singular (i.e. not stable) point. Steps in this direction have been made only in some very special cases, already quite difficult, by Laszlo \cite{laszlo} and Serman \cite{serman}. Thus the following is an important open problem about the structure of the moduli space.

\begin{problem}
Find a useful description of the tangent cone to $U_X(r,d)$ at a singular point.
\end{problem}

\subsection{Verlinde bundles and bounds for $U_X (r,d)$}

Essentially the same effective basepoint-freeness and very ampleness bounds hold replacing $SU_X(r,d)$ by $U_X(r,d)$. This however does not follow immediately, since this time one has to worry about the determinant of the vector bundles corresponding to the generalized theta divisors used in the proof. Instead, one applies results specific to the study of sheaves on abelian varieties to some natural vector bundles on the Jacobian of $X$. These were introduced in \cite{popa3} under the name of \emph{Verlinde bundles}, and are further studied in \cite{oprea1} and \cite{oprea2}.

Specifically, let $M$ be a line bundle of degree $d$ on $X$ and $\pi_M$ the composition
$$\pi_M :U_X(r,d)\overset{{\rm det}}{\longrightarrow} {\rm Pic}^d(X)
\overset{\otimes M^{-1}}{\longrightarrow} J(X).$$
Set for simplicity $U:=U_X(r,d)$ and $J:=J(X)$. Fix $F \in U_X(r_0, r_0(g-1)-d_0)$ which admits a generalized theta divisor; thus $\Theta_F$ is a \emph{basic} theta divisor. Put
$$E_m:={\pi_M}_{*}\OO_U(m\Theta_F),$$
where I simplify the notation, but remember the choice of $r,d, M, F$.

The fibers of $\pi_M$ are the $SU_X(r,L)$, for $L$ running over all line 
bundles of degree $d$ on $X$, and $\OO_U(\Theta_F)_{|SU_X(r,L)}=\LL$, the 
determinant line bundle on $SU_X(r,L)$. Kodaira vanishing (for varieties with rational singularities) and base change imply that $E_m$ is a vector bundle of rank equal to the Verlinde number $s_{r,d,m}$.
One of the main properties of the Verlinde bundles is the following: if $\phi_r: J(X) \rightarrow J(X)$ denotes the multiplication by $r$ map on the Jacobian of $X$, then 
\begin{equation}\label{decomposition}
\phi_r^* E_m \cong \underset{s_{r,d,k}}{\bigoplus}\OO_{J}(mrr_0 \Theta_{N}),
\end{equation}
where $N\in {\rm Pic}^{g-1}(X)$ is a line bundle such that $N^{\otimes r}\cong L\otimes 
({\rm det}F)^{\otimes h}$ (this follows from Exercise \ref{pullback} and the push-pull formula). Such bundles are called \emph{semihomogeneous}; an alternative description of this property is that for 
every $x\in J(X)$ there exists $\alpha \in {\rm Pic}^0 (J(X))$ such that $t_x^* E_m \cong E_m \otimes 
\alpha$. Equation (\ref{decomposition}) implies that the $E_m$ are well understood cohomologically. 

On the other hand, on an abelian variety $A$ there are precise cohomological criteria for the global generation of coherent sheaves or the surjectivity of multiplication maps, relying on the use of the Fourier-Mukai functor 
$$\R \widehat \SS: \D(A) \rightarrow \D(\widehat A), ~~\R \widehat \SS (\F) := 
\R {p_{\widehat A}}_* (p_A^* \F \otimes P),$$ 
where $P$ is a Poincar\'e bundle on $A \times \widehat 
A$. A coherent sheaf $\F$ on $X$ satisfies $IT_0$ (the Index Theorem with index $0$) if $\R \widehat \SS (\F) =  R^0 \widehat \SS (\F)$, which is equivalent to 
$H^{i}(\F\otimes\alpha)=0$ for all $\alpha\in {\rm Pic}^{0}(A)$ and all $i>0$. From our point of view here,  the key result is the following -- it was  proved first for locally free sheaves in \cite{pareschi}, and 
then for arbitrary coherent sheaves in the more general context of $M$-regularity in \cite{pp1}. 

\begin{theorem}\label{M_reg}
Let $\F$ be a coherent sheaf on $A$ satisfying $IT_0$. Then, for any ample line bundle $L$ on $A$, 
the sheaf $\F \otimes L$ is globally generated.
\end{theorem}

In the context of Verlinde bundles, the decomposition formula $(\ref{decomposition})$ allows one to  
check the $IT_0$ property for the sheaves
\begin{itemize}
\item $E_m \otimes \OO_J(- \Theta_N)$, for $m \ge h+1$. 
\item $E_m \otimes \II_{a} \otimes \OO_J(-\Theta_N)$, where $a \in J(X)$ is any point, for $m\ge 2h+1$.
\item $(E_m~ \widehat{*} ~E_m) \otimes \OO_J(-\Theta_N)$, where $ \widehat{*}$ denotes the skew-Pontrjagin product  (cf. \cite{pareschi}), for $m \ge 2h +1$.
\end{itemize}
In combination with Theorem \ref{M_reg}, this implies:

\begin{proposition}(\cite{popa3}, \cite{ep})\label{regularity}
(1)  $E_m$ is globally generated if and only if $m \ge h+1$.

\noindent
(2) $E_m \otimes \II_a$ is globally generated for $m \ge 2h +1$, for any point $a \in J(X)$.

\noindent
(3) For $m \ge 2h+1$, the multiplication map 
$$H^0 (E_m) \otimes H^0 (E_m) \longrightarrow H^0 (E_m^{\otimes 2})$$ 
is surjective. 
\end{proposition}

\begin{remark} 
In \cite{popa3}, I incorrectly claimed that for $m = r\cdot p$, the Verlinde bundles $E_m$ 
decompose as direct sums of $\OO(p\Theta)$.\footnote{Note that this also led to an erroneous statement in \emph{loc.cit.}, Theorem 5.3, namely the Claim that $\OO_U(r\Theta)$ is \emph{not} globally generated on $U_X(r,0)$; cf. also Remark \ref{dragos}. I thank D. Oprea for pointing this out.} In fact, in \cite{oprea1} Oprea shows by an explicit calculation that on elliptic curves, as part of a general 
result, for $r|m$ a decomposition into a direct sum of line bundles does exist, but the factors are various twists of $\OO(p\Theta)$ by different torsion line bundles. More recently, he has obtained for any genus  
and $d= 0$ a decomposition of $E_m$ into simple semihomogeneous factors, at least when 
$h = (r, m)$ is odd. This is discussed in Remark \ref{verlinde_other} below.
\end{remark}

Based on Proposition \ref{regularity}  and the results for $SU_X(r,d)$ in Theorems \ref{bpf} and \ref{very_ampleness}, one can then deduce via standard fibration arguments the corresponding results for $U_X(r,d)$, together with projective normality statements (see \cite{popa3} \S5,6) which I will not include here.

\begin{theorem}
Let $\Theta_F$ be a basic generalized theta divisor on $U_X(r,d)$.

\noindent
(1) (\cite{popa3}) The linear series $|m\Theta_F|$ is basepoint-free
for $m\geq {\rm  max}\{[\frac{(r+1)^{2}}{4r}] h,[\frac{r^{2}}{4}]\frac{h}{s}\}$. 

\noindent
(2) (\cite{pr}) Slightly better, on $U_X(r,0)$ the linear series $|m\Theta_F|$ is basepoint-free  for 
$m\geq {\rm  max}\{[\frac{r^{2}}{4}], r+1\}$. 

\noindent
(3) (\cite{ep}) For each $m\geq r^2+r$, the linear series 
$|m\Theta_F|$ on $U_X(r,d)$ separates points, and separates tangent vectors on the smooth locus $U^s_X(r,d)$. In particular, if $r$ and $d$ are coprime, then $|m\Theta_F|$ is very ample.
\end{theorem}

\begin{remark}\label{dragos}
The following fact, which was communicated to me by D. Oprea, shows that one is at least sometimes able to do better then the bound $r+1$ in (2) above, on $U_X(r,0)$. 

\begin{proposition}(D. Oprea)
If $|\LL|$ is basepoint-free on $SU_X(r)$, then $|2\Theta_L|$ is basepoint-free on $U_X(r,0)$ for any $L\in {\rm Pic}^{g-1}(X)$. In particular, by Theorem \ref{morphism}, this holds when $r=2$ and when $r=3$ and $X$ is generic.
\end{proposition}
\begin{proof}
This goes essentially as in the case of abelian varieties. First, by the geometric interpretation of Strange Duality, Remark \ref{interpretation}, $|\LL|$ is generated by theta divisors, which can be easily translated as saying that for each $E$ in $U_X(r,0)$, 
$H^0 (E \otimes M) = 0$ for $M \in {\rm Pic}^{g-1}(X)$ general. Now if $L \in {\rm Pic}^{g-1}(X)$ is any line bundle, then in the linear series $|2\Theta_L|$ on $U_X(r,0)$ we have 
all reducible divisors of the form $\Theta_{L\otimes A} + \Theta_{L\otimes A^{-1}}$, for
$A \in {\rm Pic}^0 (X)$. By the observation above, any fixed $E$ is avoided by such a 
divisor if $A$ is general.
\end{proof}
\end{remark}

The same techniques, combined with the sharp results on very ampleness for rank $2$ bundles with fixed determinant mentioned at the beginning of \S7, imply in rank $2$ a better result.

\begin{theorem}(\cite{ep})
Let $\Theta_F$ be a basic generalized theta divisor on $U_X(2,d)$. If $d$ is odd, then $|3\Theta_F|$ is very ample. If $d$ is even, and $C$ is not hyperelliptic,  then $|5\Theta_F|$ is very ample.
\end{theorem}

\begin{remark}[{\bf Decomposition, and other uses of the Verlinde bundles}; cf. \cite{oprea2}, \cite{popa3}]\label{verlinde_other}
(1) The Fourier-Mukai transform of the Verlinde bundles is an equally useful object in 
terms of linear series on moduli spaces. For simplicity, let's assume here that $d=0$, and redenote $E_m$ by $E_{r,m}$, to include the rank. The Strange Duality isomorphism (cf. \S5.2), in fact their collection over all moduli spaces $SU_X(r,L)$, can be interpreted as a natural isomorphism
\begin{equation}\label{global_sd}
SD: E_{r,m}^\vee \longrightarrow \widehat{E_{m,r}},
\end{equation}
where the hat denotes the Fourier-Mukai transform coming from the dual Jacobian. One can for instance show that for $m=1$ this is an isomorphism by checking the stability of the bundles involved
(see \cite{popa3} \S3). Furthermore, Oprea shows in \cite{oprea2} the following statement on the decomposition of  $E_{r,m}$ into simple semihomogeneous (hence stable by \cite{mukai1}) factors, extending his previous result for elliptic curves in \cite{oprea1}.

\begin{theorem}[\cite{oprea2}, Theorem 1]\label{decomposition2}
If $(r, k) =1$ and $h$ is odd, the Verlinde bundle $E_{hr, hk}$ splits as
$$E_{hr, hk} \cong \underset{\xi}{\bigoplus} W_{r,k, \xi}^{\oplus m_{\xi} (r,k)}$$ 
where $\xi$ runs over all $h$-torsion line bundles on $J(X)$ .
\end{theorem}

\noindent
The notation is as follows:  Oprea deduces from results in Mukai \cite{mukai1} that on a principally polarized abelian variety $(A, \Theta)$ with $\Theta$ symmetric, for every coprime integers $r$ and $k$ there exists a unique symmetric simple semihomogeneous bundle $W_{r,k}$ such  that
$${\rm rk}~W_{r,k} = r^g, ~~~~{\rm det} ~W_{r,k} = \OO_A (r^{g-1}k \Theta).$$
(Symmetric means $(-1)^* W_{r,k} \cong W_{r,k}$; without this requirement uniqueness is not 
satisfied.) One can easily check that for every $\xi\in {\rm Pic}^0 (J(X))$ and $\bar{\xi}$ an $r$-th
root of $\xi$, the bundle $W_{r,k} \otimes \bar{\xi}$ depends only on $\xi$. We then denote 
$$W_{r,k, \xi} := W_{r,k} \otimes \bar{\xi}.$$
The multiplicities $m_{\xi}(r, k)$ are computed explicitly in \cite{oprea2}, by means of a calculation of 
the trace of the action of $h$-torsion points on the spaces of generalized theta functions (itself very interesting), but I will not reproduce the complicated formula here. They depend 
on the Verlinde numbers on a higher genus curve (see \emph{loc. cit.}, Theorem 2).
One can check that $W_{r,k,\xi}^\vee \cong \widehat{W_{k,r,\xi}}$ and 
$m_{\xi} (r,k) = m_{\xi} (k,r)$, hence the decomposition in the Theorem agrees term by term 
with the Strange Duality isomorphism (\ref{global_sd}). This does not however reprove Theorem \ref{sd}, as the various isomorphisms are not known to be compatible with the map in (\ref{global_sd}).

\noindent
(2) Finally, with suitable choices of $r$, $m$, and $g$, it is shown in \cite{popa3} that the restriction of 
$\widehat{E_{m,r}}$ to a suitable Abel-Jacobi embedding of $X$ provides a class of base points for the series $|\LL|$ on $SU_X(r)$ extending Raynaud's examples discussed in \S6.2. However, 
as the decomposition in Theorem \ref{decomposition2} is now available, let me note that this is really 
a statement about the bundles $W_{r,k}$ (and their twists). Indeed, precisely as in \cite{popa3} \S4, 
using (16) in \cite{oprea2} instead of (10) here, one sees that the restriction of $W_{r,m}^\vee$ 
to a suitable embedding of $X$ in its Jacobian is a semistable vector bundle of rank $m^g$ and 
slope $\frac{kg}{m}$, satisfying condition $(*)$ in \S6.2 when $k <m$ and $m|g$. Raynaud's examples 
are exactly the case $k = 1$; note that the argument works for the Fourier-Mukai transform of 
every simple semihomogeneous bundle.
\end{remark}

\subsection{Conjectures}

It is likely that various effective bounds mentioned above are not optimal. In particular one hopes for linear bounds in the rank. Inspired by the low-rank cases and the uniform bound for the global generation of Verlinde bundles in \ref{regularity}, I would
speculate the following: 

\begin{conjecture}\label{bounds}
(a) If $\Theta$ is a basic theta divisor, the linear series $|m \Theta|$ on $U_X(r,0)$ is 
\begin{enumerate}
\item basepoint-free for $m \ge r+1$.
\item very ample for $m \ge 2r + 1$.
\end{enumerate}

\noindent
(b) The linear series $|\LL^m|$ on $SU_X(r)$ is
\begin{enumerate}
\item basepoint-free for $m \ge r-1$.
\item very ample for $m \ge 2r + 1$.
\end{enumerate}

\noindent
(c) The bounds are at least as good on $U_X(r,d)$ and $SU_X(r,d)$ with $d \neq 0$.
\end{conjecture}

\begin{remark}
(1) In the Conjecture, (a) trivially implies (b), with $r+1$ instead of $r-1$ in the case of basepoint-freeness, by restriction. 

\noindent
(2) More importantly, (b) implies (a) via the same procedure using Verlinde bundles on Jacobians.

\noindent
(3) The main appeal of the conjecture in (a) is that it directly extends the case of classical theta functions of Jacobians (and more generally ample divisors on any abelian variety) given by the Lefschetz theorem, as a linear function in the rank. As we see for example in Proposition \ref{dragos}, this 
would still not always be optimal, but at the moment no general machinery seems to indicate towards a better guess. The conjecture in (b) is a bit more speculative. Also, there is again the possibility that one can do even better. For instance, Beauville conjectures in \cite{beauville5} that $|\LL|$ is basepoint-free on $SU_X(3)$, which would be below this bound. It is even harder to predict a bound for very ampleness on $SU_X(r)$, so somewhat inconsistently I haven't even attempted to go beyond what is implied by (a), namely that it should be better than $2r+1$. In the same vein, it is not yet clear what sharper conjecture should replace (c) above in the case of arbitrary degree.
\end{remark}

\end{document}